\documentclass[10pt]{article}
\usepackage[latin1]{inputenc}
\usepackage{graphicx}
\usepackage{amssymb}
\usepackage[all]{xy}
\usepackage{amsfonts}
\usepackage[thmmarks]{ntheorem}

\def\dcg{[\![}
\def\dcd{]\!]}

\def\H{{\mathcal H}}

\def\C{{\mathcal C}}

\def\S{{\mathcal S}}

\def\RR{{\mathbb R}}
\def\CC{{\mathbb C}}
\def\NN{{\mathbb N}}
\def\LL{{\mathcal L}}
\def\interior{{\rm Int}}
\def\reel{{\rm Re}}

\def\bra{\langle}
\def\ket{\rangle}
\def\id{{\rm Id}}
\def\tr{\mbox{\rm Tr}}

\def\bea{\begin{eqnarray}}
\def\eea{\end{eqnarray}}

\def\be{\begin{equation}}
\def\ee{\end{equation}}

\def\Herm{{\rm Herm}}
\def\diag{{\rm diag}}

\newtheorem{theorem}{Theorem}
\newtheorem{definition}{Definition}
\newtheorem{lemma}{Lemma}
 
\newtheorem{propo}{Proposition}
\newtheorem{cor}{Corollary}
\theoremstyle{nonumberplain}
\theorembodyfont{\normalfont}
\theoremseparator{:}
\theoremsymbol{$\P$}
\newtheorem{demo}{Proof}

\newcommand{\bigL}{\mathop{\mathrm{Lex}}}
\begin{document}
 
% \begin{frontmatter}

\title{Noncommutative ordered spaces : examples and counterexamples}
\author{Fabien Besnard\footnote{Pôle de recherche M.L. Paris, EPF, 3 bis rue Lakanal, F-92330 Sceaux. \goodbreak fabien.besnard@epf.fr}}
% \address{Pôle de recherche M.L. Paris, EPF, 3 bis rue Lakanal, F-92330 Sceaux\\ fabien.besnard@epf.fr}

\maketitle
\begin{abstract}
In order to introduce the notion of causality in noncommutative geometry it is necessary to extend Gelfand theory to the context of ordered spaces. In a previous work we have already given an algebraic caracterization of the set of non-decreasing continuous functions on an certain class of topological ordered spaces. Such a set is called an isocone, and there exist at least two versions of them (strong and weak) which coincide in the commutative case. In this paper we introduce yet another breed of isocones, ultraweak isocones, which has a simpler definition with a clear physical meaning. We show that ultraweak and weak isocones are in fact the same, and completely classify those which live in a finite dimensional $C^*$-algebra, hence corresponding to finite noncommutative ordered spaces. We also give some examples in infinite dimension.
\end{abstract}
\section{Introduction}
The application of noncommutative geometry to particle physics initiated by Connes and Chamseddine in the 90's is one of the most promising ideas to shed some light on the seemingly arbitrary ingredients of the standard model (see \cite{ChCoMa} for a detailed explanation, and \cite{resil}, \cite{ChCoVS1}, \cite{ChCoVS2} and \cite{DLM} for some recent and important updates on the subject). It can loosely be described as a ``noncommutative Kaluza-Klein'' model, where the geometry of spacetime splits into a direct product of an ordinary manifold $M$ and a metaphoric noncommutative finite space $F$, a so-called ``almost commutative geometry''. A purely geometric action on $M\times F$ then yields the complete bosonic part of the standard model action minimally coupled to gravity, hence realizing the unification of all forces in a framework which is very close in spirit to general relativity. However, it was clear from the beginning that there are some inherent limitations to this model. First, it does not include any quantum gravity effect. In other words, the manifold $M$ remains commutative and smooth at all energy scales. This is an important step, but certainly a preliminary one, and on this matter we would like to quote \cite{spectralpov} (emphases are in the original text):

\begin{quote}
``It could be well that the coherence of the spectral action principle indicates that our continuum picture of space-time is only an approximation to a completely \emph{finite} spectral geometry whose underlying Hilbert space is \emph{finite dimensional}. [\ldots] In this scenario, once we go up in energy towards the unification scale, the small amount of noncommutativity encoded in the finite geometry $F$ to model the present scale, will gradually creep in and invade the whole algebra at Planck scale.'' (A. H. Chamseddine, A. Connes, 2010)
\end{quote}

The second limitation comes from the signature of the metric, which is only allowed to be euclidean in the current formulation of the model. In this case talking about spacetime as we have done above is a bit of an overstatement. The problem can be circumvented by a procedure known as Wick rotation, but only when the curvature of the smooth manifold vanishes, which can always be assumed in current experimental particle physics. However, the very fact that the physical signature of the metric must be Lorentzian has important repercussions in the finite geometry as well (\cite{barrett}). Clearly something has to be done about this issue.

There already exist several approaches to ``Lorentzian noncommutative geometry''. We recommend \cite{franco} for a thorough review. The path we have chosen is to focus on causality. Indeed, the Lorentzian signature is singled out among every others by the fact that it allows for the definition of a partial order structure on the set of events, at least locally. Moreover, knowledge of the causal structure permits to recover the metric up to a conformal factor, which is a local scale for the measurement of durations. That these two aspects of time, duration and causality, play quite distinct roles is an important lesson of relativity, it thus  seems natural to split the degrees of freedom in this way. As a matter of fact, it is in essence the point of departure of the causal set approach to quantum gravity (\cite{causet}). We therefore propose the following sketch of a program :  define causality in noncommutative geometry, incorporate the conformal factor, write down the dynamical equations of the theory. Let it be clear that in this paper we will only deal with the first step.

From a purely mathematical perspective, that we adopt in most of what follows, causality is just a partial order relation. In view of the quotation above, what we have to do is clear, if not straightforward : define what a noncommutative partially ordered space should be, and investigate with particular care the case of finite spaces (i.e. finite-dimensional algebras).  In \cite{bes1} we have already given a first tentative definition of a partial order structure on a noncommutative space. In fact, we have given two of them. Moreover, we have given a class of examples within the algebra $M_2(\CC)$, where the two definitions happen to agree. However several issues remained. First, we had two definitions instead of one. Moreover, neither of them was very palatable, nor easy to work with. Finally we needed to find more examples to show the interest of the definition(s). In this paper we will address all these issues.

First, we will pick one of the two definitions as our preferred choice by showing that it is in fact equivalent to a third, much more natural and tractable. Then we will show how to build new noncommutative ordered spaces out of already known ones. One of these constructions will be to add ``noncommutative infinitesimals''. This will provide us with  a completely different family of examples from the one studied in \cite{bes1}. Finally we will state and prove a classification theorem in the finite-dimensional case.

The paper is organized as follows : in section 2 we recall the  definitions of strong and weak isocones and $I^*$-algebras, and summarize the results already obtained in \cite{bes1}. In section 3 we introduce a new class of isocones : ultraweak isocones, and show its equivalence with weak isocones. From this point on, we only consider weak isocones. Section 4 is devoted to some constructions involving isocones and several examples. Finally we study finite-dimensional $I^*$-algebras in section 5 and give their classification.

Here are some notations which we will use throughout the text. The $C^*$-algebra of continuous complex valued functions defined on the compact set $M$ will be written as ${\cal C}(M)$. If $A$ is a $C^*$-algebra we denote by $\reel(A)$ the set of self-adjoint elements of $A$, and by $A^+$ the set of its positive elements. This rule admits the following two exceptions : the real part of ${\cal C}(M)$ will be written ${\cal C}(M,\RR)$, and the set of hermitian $N\times N$ matrices will be denoted by $\Herm(N)$.

All our $C^*$-algebras have a unit, and sub-$C^*$-algebras contain the unit. The spectrum of $a$ is written $\sigma(a)$. If $H$ is a Hilbert space, ${\cal B}(H)$ will be the $C^*$-algebra of bounded operators on $H$, and ${\cal K}(H)$ its ideal of compact operators.

We will use freely and frequently the following facts about the continuous functional calculus : it commutes with $*$-morphisms, and is continuous in its operator argument, i.e. if $f\in{\cal C}(\RR)$, the map $a\mapsto f(a)$ defined on $\reel(A)$ is continuous (\cite{halmos}, problem 126). 

Here are some notations concerning the matrix algebra $M_N(\CC)$ and the direct sum $S=\bigoplus_{x=1}^kM_{n_x}(\CC)$. The unit of $M_N(\CC)$ will be written $1_N$, and its zero $0_N$. The projection of $S$ onto its $x$-th summand is $\pi_x$, elements of $S$ are written $(a_x)_{1\le x\le k}$. The element of $S$ such that $a_y=1_{n_y}$ if $y=x$ and else $a_y=0_{n_y}$ will be written $\iota_x$.  We recall that a hermitian matrix which has at least one multiple eigenvalue is said to be derogatory. The manifold of non-derogatory matrices in $\Herm(N)$ will be written $S_N$. When needed we will use the Frobenius inner product $\bra a,b\ket:=\tr(ab^*)=\tr(ab)$ for $a,b\in \Herm(N)$.

\section{Isocones and $I^*$-algebras}
In noncommutative geometry, one trades spaces for algebras. What is meant by ``spaces'' varies, but at the very least these are Hausdorff locally compact topological spaces. In this paper we will assume all spaces to be compact, for simplicity's sake\footnote{Some might worry that this would exclude causal Lorentz manifolds. However causal Lorentz manifolds can be compact provided they have a boundary. At the topological level, which is the one we work in,  adding a boundary causes no problem.  We refer to \cite{bes1} about the compactifications procedures in the topological ordered space setting.}. By ``algebras'' we mean commutative unital $C^*$-algebras, and by ``trade'' we mean that the category of compact Hausdorff spaces equipped with their continuous mappings is dually equivalent to the one of commutative unital $C^*$-algebras with their $*$-morphisms. This is the content of the well-known (commutative) Gelfand-Naimark theorem. Once we remove the adjective ``commutative'', we cannot trade the algebras back for spaces, and we are entering the realm of noncommutative geometry (the word ``topology'' would be more appropriate unless we have more structure).

Starting with a compact Hausdorff space $M$, we now want to add a partial ordering $\preceq$ on it. Of course the topology and the partial order have to be compatible in a way. The functions which are sensitive to both structures are those which are continuous and order-preserving, hence we introduce some names and notations for them.

\begin{definition} Let $M$ (resp. $N$) be a topological space equipped with a partial order $\preceq$ (resp. $\le$). A map $f$ from $M$ to $N$ is \emph{isotone} iff it satisfies 

\be
x\preceq y\Rightarrow f(x)\le f(y)
\ee
 for all $x,y\in M$.  We denote the set of maps from $M$ to $N$ which are both continuous and isotone by $I(M,N)$. When $N$ is $\RR$ equipped with the natural ordering, we   put $I(M):=I(M,\RR)$.
\end{definition}

The elements of $I(M)$ will simply be called  real isotonies, the continuity being understood. In order to expect a duality result of the Gelfand-Naimark sort, we have to suppose that there are ``enough'' real isotonies.

\begin{definition} A topological ordered set $M$ is said to be \emph{completely separated} iff for all $x,y\in M$, 
\be
x\preceq y\Leftrightarrow \forall \phi\in I(M),\ \phi(x)\le\phi(y)
\ee
\end{definition} 

A completely separated topological ordered space will be called a \emph{toposet} in the rest of the paper. We recall that if $M$ is compact, it is a toposet if and only if the order relation $\preceq$ is closed in $M\times M$ for the product topology.

It turns out that it is possible to characterize algebraically the sets $I(M)$ where $M$ is a compact toposet, and to recover the toposet from the algebraic structure. We need some more definitions to be more specific. Recall that in any $C^*$-algebra, the meet and join of two elements can be defined through functional calculus by the formulas :

\be
a\vee b:={1\over 2}(a+b+|a-b|),\quad a\wedge b:={1\over 2}(a+b-|a-b|)
\ee

Note that these operations do not satisfy the lattice axioms when the elements $a$ and $b$ do not commute. Of course, they reduce to the usual supremum and infimum of functions when the algebra is commutative.

\begin{definition}\label{defisoc} Let $A$ be a $C^*$-algebra with unit $1$.  A subset $I$ of $\reel (A)$ which satisfies :
\begin{enumerate}
\item $\forall x\in\RR$, $x.1\in I$.\label{constants}
\item $\forall b,b'\in I$, $b+b'\in I$,\label{somme} 
\item $\forall \lambda\in\RR_+$, $\forall b\in I$, $\lambda b\in I$,\label{scal}
\item $\forall b,b'\in I$, if $b$ and $b'$ commute then $b\vee b'\in I$ and $b\wedge b'\in I$,\label{meetandjoin}
\item $\overline{I}=I$ ($I$ is norm-closed),\label{ferme}
\end{enumerate}

will be called a (weak) pre-isocone. If \ref{meetandjoin} holds also when $b$ and $b'$ do not commute, it is called a strong pre-isocone. A weak (resp. strong) pre-isocone will be called a weak (resp. strong) isocone if it moreover satisfies 
\begin{enumerate}\setcounter{enumi}{5}
\item\label{dense} $\overline{I-I}=\reel(A)$
\end{enumerate}

A couple $(I,A)$ when $I$ is a weak (resp. strong) isocone of $A$ is called a weak (resp. strong) $I^*$-algebra\footnote{The reason for this name is the following : the $C$ of $C^*$-algebra stands for ``continuous'', while the $I$ of $I^*$-algebra stands for ``isotony''.}.
\end{definition}

Note that the closure of the set of differences of elements of $I$ appearing in axiom \ref{dense} is the same as the closure of the linear span of $I$, since $I$ is convex. 

By default, an isocone (pre-isocone, $I^*$-algebra) will be considered to be of the weak sort.  An obvious example of strong isocone is the \emph{trivial isocone} $I=\reel (A)$.

The following theorem justifies the above definitions :

\begin{theorem} The Gelfand transform realizes a dual equivalence of categories between the category of commutative $I^*$-algebras with their morphisms and the category of compact toposets with their continuous isotonies.
\end{theorem}

We refer to \cite{bes1} for the proof of this result and all the others in this section. Let us just say that $I^*$-morphisms are $*$-morphisms mapping the isocone of one algebra into the isocone of another, and satisfying some extra-conditions that will not concern us here (and which are trivially satisfied in the commutative case). We also recall that the Gelfand transform associates to a commutative $C^*$-algebra $A$ its space of characters ${\cal X}(A):=\{\phi : A\rightarrow \CC|\phi$ is a $*$-morphism$\}$, and to an element $a\in A$ the function $\hat a : {\cal X}(A)\rightarrow \CC$, such that $\hat a (\phi)=\phi(a)$ for every $\phi\in{\cal X}(A)$. The partial order structure on ${\cal X}(A)$ is naturally defined by the isocone $I$ in $A$ through 

\be 
\phi\preceq_I\psi\Leftrightarrow \forall a\in I,\ \phi(a)\le\psi(a)\label{stateordering}
\ee

Let us   observe that the above formula   defines a toposet structure on some of the various spaces attached to $A$   when $A$ is noncommutative.

\begin{propo} Let $(I,A)$ be an $I^*$-algebra. The formula (\ref{stateordering}) defines a toposet structure on :  the set of states $S(A)$, the set of pure states $P(A)$ and the set of characters ${\cal X}(A)$. 
\end{propo}

In fact $\preceq_I$ is none other than the ordering on the hermitian forms on  $A$ defined by the dual cone $I^*=\{\omega\in \reel(A^*)|\omega(I)\subset\RR_+\}$, where $\reel(A^*)$ denotes the set of hermitian forms on $A$. Since $1_A$ and $-1_A$ belong to $I$, we see that $I^*$ lies inside the hyperplane defined by the equation $\omega(1_A)=0$. Hence, hermitian forms are comparable with respect to $\preceq_I$ only if they have the same value on $1_A$, which is of course the case for the states on $A$.

Another noteworthy result is the following.

\begin{theorem}\label{isotcalc}
Let $(I,A)$ be a strong or weak $I^*$-algebra, and $a\in I$. Then $I\cap C^*(a)$ is an isocone in $C^*(a)$. Under the identification of ${\cal X}(C^*(a))$ with  $\sigma(a)$, $I\cap C^*(a)$ induces through (\ref{stateordering}) a toposet structure $\preceq_a$ on $\sigma(a)$ which is at most as fine as the natural ordering of $\RR$. Therefore $I$ is ``stable by isotone calculus'', that is to say :
\be
\forall a\in I,\  \forall f\in I(\sigma(a)),\ f(a)\in I\label{isotonecalc}
\ee 
\end{theorem}

Here $I(\sigma(a))=I(\sigma(a),\le)$ where $\le$ is the natural order of $\RR$. Of course the stability under isotone calculus is true also for $f\in I(\sigma(a),\preceq_a)$ (and $I(\sigma(a))\subset I(\sigma,\preceq_a)$).

Observe that since every continuous increasing function on the compact space $\sigma(a)$ can be extended to a continous increasing function on $\RR$, the set $I(\sigma(a))$ can as well be replaced by the set of real isotonies $I(\RR)$ in (\ref{isotonecalc}). Hence an isocone (weak or strong) satisfies the stability property :

\be
\forall f\in I(\RR),\  f(I)\subset I
\ee

The isotone functional calculus has a consequence which will prove to be particularly important in the  finite-dimensional case.

\begin{propo}\label{combproj} Let $I$ be an isocone in the finite-dimensional $C^*$-algebra $A$. Then for all $a\in I$, there exist commuting projections $p_1,\ldots,p_n\in I$, positive constants $\lambda_1,\ldots,\lambda_n$ and some real constant $\lambda$ such that 
$$a=\sum_{i=1}^n\lambda_ip_i+\lambda.1_A$$
\end{propo}

Before turning to examples, we recall another corollary of isotone functional calculus that will be used in this paper. 

\begin{propo} Let $I$ be an isocone, and $c_1,\ldots,c_n$ be a family of positive operators in $I$ commuting among themselves. Then the product $c_1\ldots c_n\in I$.
\end{propo}

Finally we recall a family of examples of isocones in the simplest noncommutative $C^*$-algebra, namely $M_2(\CC)$.

In the case $N=2$, it turns out that any closed convex cone containing $\RR.1_2$ is automatically stable by noncommutative $\vee$ and $\wedge$, and is therefore a strong pre-isocone. We can thus classify isocones with $M_2(\CC)$ as algebra by coordinatizing such cones with non-empty interior. Here is a convenient way to do it. Let $S$ be the sphere of radius $1$ and center $1_2/2$ for the Frobenius norm inside  $H_1:=\{a\in\Herm(2)|\tr(a)=1\}$ (hence $S$ is a $2$-sphere). We note that $S$ is also the set of rank $1$ projectors in $\Herm(2)$. For any compact subset  $K$ of $S$ we set $I_K:=\RR_+.K+\RR.1_2$. We say that $K$ is \emph{geodesically convex} iff for any two non-antipodal points $x,y\in K$, the smallest arc of great circle joining $x$ to $y$ lies in $K$. We then have the following proposition.

\begin{propo} The isocones in $M_2(\CC)$ are the sets $I_K$ where $K$ is either $S$ or a compact geodesically convex subset of $S$ with non-empty interior (for the relative topology of $S$).
\end{propo}

It is also possible to prove that $I^*$-isomorphisms correspond to isometries of $S$. Let us end by a simple illustrative example of the previous theorems. Take $a\in I_K$ with two distinct eigenvalues $a_1<a_2$. Let us ask what is the ordering $\preceq_a$ induced on $\sigma(a)$ by $I$ ? Let $p_2$ be the projection on the eigenspace of $a_2$. Then $p_2\in I_K$ by isotone calculus. Hence  $p_2\in K$ since $K$ is the set of rank $1$ projections in $I$. Then :

\begin{itemize}
\item If $-p_2$ (or equivalently $p_1=1_2-p_2$) also belongs to $K$, the ordering $\preceq_a$ is trivial.
\item If not, then it is the natural ordering.
\end{itemize}

\begin{figure}[hbtp]
\begin{center}
\includegraphics[scale=0.7]{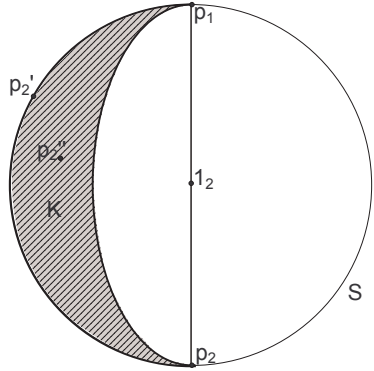}
\end{center}
\caption{In this example, if the eigenprojector on the largest eigenvalue of $a$ is $p_2$, then $\preceq_a$ is trivial. If it is $p_2'$ or $p_2''$, it is the natural order on $\sigma(a)$.}
\end{figure}

Note also that thanks to the Frobenius scalar product, the state space and the ordering $\preceq_I$ induced on it by $I$ can be ``internalized'' (the same goes of course for any finite-dimensional $C^*$-algebra). Any state $\phi$ on $M_2(\CC)$ can be written thanks to a density matrix $\rho$, that is to say a positive element of $H_1$, through  $\phi(a)=\tr(\rho a)$ for any $a\in M_2(\CC)$. Using this identification, the order relation $\preceq_I$ on density matrices can be defined by $\rho\preceq_I\rho'\Leftrightarrow \rho'-\rho\in I^*$, where $I^*$ is the internal dual cone of $I$, namely the set of hermitian (and trace-free) matrices $m$ such that $\tr(am)=0$ for all $a\in I$.

\section{Ultraweak isocones}

The duality theorem encourages us to think of general $I^*$-algebras as the noncommutative counterparts of toposets. However, there exists an infinity of noncommutative generalizations of the same commutative structure, and weak and strong $I^*$-algebras already provide two examples. In this section we will show how a ``naturalness'' criterion allows to pick one of them.  Let us begin with a definition.

\begin{definition}\label{defisoc2} Let $A$ be a $C^*$-algebra with unit $1$.  A subset $I$ of $\reel (A)$ which satisfies :
\begin{enumerate}
\item $\forall x\in\RR$, $x.1\in I$.\label{constants2}
\item $\forall b,b'\in I$, $b+b'\in I$,\label{somme2} 
\item $\forall \phi\in I(\RR),\ \phi(I)\subset I$,\label{isotstabl}
\item $\overline{I}=I$,\label{ferme2}
\end{enumerate}

will be called an ultraweak pre-isocone. If $I$ moreover satisfies 
\begin{enumerate}\setcounter{enumi}{4}
\item\label{dense2} $\overline{I-I}=\reel(A)$.
\end{enumerate}
it will be called an ultraweak isocone.
\end{definition}

Ultraweak $I^*$-algebras are defined accordingly. Let us discuss this set of axioms.

We note first that since non-decreasing linear functions belong to $I(\RR)$, ultraweak pre-isocones are necessarily convex cones, and axiom \ref{constants2} could be replaced by the requirement that $I$ is non-empty, which is redundant in case axiom 5 is satisfied.

From the physical point of view,  axiom \ref{isotstabl} is almost a tautology. Indeed, if we interpret an  element a of $I$ as some  causal observable, then $\phi(a)$ will represent the observation of $a$ followed by a ``nonlinear rescaling'' of the real line by $\phi$, an operation that has no physical consequence as far as the causal relations only are concerned. Thus $\phi(a)$ must also be a causal observable. 

Axiom \ref{ferme2} seems more like a mathematical convenience than a physically important fact. Indeed, if ever a set of causal observables $I$ satisfied all the other axioms, its norm closure would also.

Axiom \ref{dense2} is required in order to define a partial order relation on the state space, instead of just a partial pre-order.

In the end, the only axiom for which we do not find a straigthforward justification is the second one. Incidentally, we do not find more justification to the widely accepted fact in quantum physics that the sum of two observables has to be an observable also (of course we implicitly accepted this by using the $C^*$-formalism). On this issue we will take a pragmatic approach by observing that at the very least, we would have to require the set of causal observables to be stable under sum when the two terms of the sum commute. But if we required only this, then the whole theory would become trivial : the set $I$ of causal observables would just be a union $I=\bigcup I_M$ over some family of toposets $M$, with $I_M$ isomorphic to $I(M)$. Since we do not want a trivial theory, we accept the only axiom which allows noncommuting causal observables to ``interfere'' with one another, namely axiom \ref{somme2}.

\underline{Remark} : We conjecture that this set of axioms is equivalent to an even simpler and more natural one, namely that $I$ defines a partial order on the pure state space through (\ref{stateordering}), that it is stable by isotone functional calculus, and that the elements of $A$ which have a non-decreasing Gelfand transform with respect to $\preceq_I$ are all in $I$. This conjecture, which is under investigation, seems to have close ties with the Stone-Weïerstrass conjecture for $C^*$-algebras.

So these axioms are well and good since they are in some sense minimal. It would appear then that we would have to prove a duality theorem using these axioms only, and declare our older work obsolete. However we will not need to do so, since it turns out that ultraweak isocones and weak isocones are in fact one and the same thing ! Of course we know that weak isocones are ultraweak isocones by theorem \ref{isotcalc}. To prove the converse we will use this result   : 

\begin{theorem}\label{approx}
Let $(M,\preceq )$ be a compact Hausdorff partially ordered set. Let $A$ be the set of piecewise linear elements of $I(\RR)$. Let $S$ be a non empty subset of $I(M)$. If
\begin{enumerate}
\item $S$ is stable by sum,
\item\label{h2} $\forall f\in S$, $\forall \phi\in A$, $\phi\circ f\in S$, 
\item\label{h3} $\forall x,y\in M$, $x\preceq y\Leftrightarrow \forall f\in S$, $f(x)\le f(y)$.
\end{enumerate}
then $S$ is dense in $I(M)$ for the uniform norm. 
\end{theorem}

We refer to  (\cite{bes2}) for the proof. Now we can show the claimed equivalence.

\begin{theorem} $(I,A)$ is an ultraweak $I^*$-algebra iff $(I,A)$ is a weak $I^*$-algebra.
\end{theorem}
\begin{demo} All we need to prove is that if $I$ is an ultraweak isocone and $a_1,a_2\in I$, with $a_1a_2=a_2a_1$, then $a_1\wedge a_2$ and $a_1\vee a_2$ belong to $I$. We know that there exists a compact subset $M\subset \RR^2$ and a $*$-isomorphism $\Psi : C^*(a_1,a_2)\rightarrow {\cal C}(M)$ such that $\Psi(a_i)=\pi_i$, the projection on the $i$-th coordinate. Let $S=\Psi(I\cap C^*(a_1,a_2))$, and define $\preceq$ on $M$ by   $\forall x,y\in M$, $x\preceq y\Leftrightarrow \forall f\in S$, $f(x)\le f(y)$. It is obviously a preorder on $M$. Moreover, if $f(x)=f(y)$ for all $f\in S$, then $\pi_i(x)=\pi_i(y)$, $i=1,2$. Thus  $x=y$ and $\preceq$ is a partial order relation on $M$. Now since a $*$-morphism commutes with functional calculus, it is clear that $S$ satisfies the hypotheses of theorem \ref{approx}. Moreover $S$ is closed since $I$ is and $\Psi$ is a $*$-isomorphism. Since $I(M)$ is closed under $\vee$ and $\wedge$, we have $a_1\vee a_2,a_1\wedge a_2\in I$. 
\end{demo}

% For the moment we just recall a family of example found in \cite{bes1} which unfortunately cannot settle the issue since both the condition 4 of definition \ref{defisoc} and its stronger version are entailed by the other axioms in this particular case. 

In view of this theorem, the words pre-isocone, isocone, and $I^*$-algebra will always refer to the ultraweak versions in the following.%, except when specified otherwise.

% \underline{Remark} : We will not use this result, but it can be easily shown that the set of functions which stabilize a non-trivial isocone $I$ is exactly $I(\RR)$.

\section{Some operations on isocones}
Let us recall some concepts and notations of order theory. For any two elements $x,y$ in a poset $(M,\preceq)$, we say that $x$ and $y$ are comparable iff $x\preceq y$ or $y\preceq x$, and we write $x\perp y$. If they are incomparable we write $x\parallel y$. If $x\preceq y$ and $x\not=y$ we write $x\prec y$ and say that $x$ is strictly below, or strictly less than $y$. We can also extend the meaning of the symbols $\preceq$, $\prec$ and $\parallel$ to subsets of $M$. For instance if $X$ and $Y$ are subsets of $M$, we write $X\prec Y$ if every element of $X$ is strictly below every element of $Y$, and $X\parallel Y$ if no element of $X$ is comparable to any element of $Y$.

Given two posets $(M,\preceq_M)$ and $(N,\preceq_N)$, we can form their \emph{disjoint sum} or \emph{cardinal sum}, denoted by $M+N$, which is equal to $M\coprod N$ as a set, with relations $\preceq_M$ on $M$, $\preceq_N$ on $N$, and $M\parallel N$. If one replaces the relation $M\parallel N$ by $M\prec N$ in the previous construction, one obtains the \emph{ordinal sum} of $M$ and $N$, which is denoted by $M\oplus N$. Of course this is not a very good notation, since this operation is not commutative, but since it is standard we will stick to it .

The operations of cardinal and ordinal sums allow us to write many\footnote{But not all: a poset whose Hasse diagram is an ``N'' cannot be so obtained.} finite posets in a convenient way. For instance we can write $(1+2)\oplus (3+4)$. This means that we consider the poset with elements $\{1;2;3;4\}$ subjects to the relations $1\preceq 3$, $1\preceq 4$, $2\preceq 3$, $2\preceq 4$, plus the obvious $x\preceq x$ for every $x$. This poset can also be visualized thanks to its nicer looking but longer to typeset Hasse diagram :

% $$\xymatrix@R-15pt{
%  & 1+t & \\
%  & S^+\ar[u]& \\
% 1\ar[uur]&  &t\ar[uul]\\
%  & S^-\ar[uu]&\\
%  & 0 \ar[uul]\ar[u]\ar[uur]&
% }
% $$ 

$$\xymatrix {
 3& & 4\\
1\ar@{-}[u] \ar@{-}[urr]& &2\ar@{-}[u] \ar@{-}[ull]
}
$$

This is the graph of the covering relation ($y$ covers $x$ iff $x\prec y$ and there is no $z$ such that $x\prec z\prec y$) with the convention that elements cover others from top to bottom. 

We now consider a poset $P$ and a family of posets $(M_x, \preceq_x)_{x\in P}$. The \emph{lexicographic sum} of the posets $M_x$ over $P$ is the disjoint union ${\displaystyle \coprod_{x\in P} M_x}$ equipped with the order relation $\preceq$ defined by : $\preceq$ is equal to $\preceq_x$ when restricted to elements of $M_x$, $x\prec y\Rightarrow M_x\prec M_y$ and $x\parallel y\Rightarrow M_x\parallel M_y$. The lexicographic sum will be denoted by ${\displaystyle \bigL_{x\in P}M_x}$. Given the posets $1+2$ and $1\oplus 2$, the lexicographic sum supersedes cardinal and ordinal sums : $M_1+M_2=\bigL_{i\in 1+2}M_i$, and $M_1\oplus M_2=\bigL_{i\in 1\oplus 2}M_i$.

We now carry  over the lexicographic sum to $I^*$-algebras.

\begin{theorem}\label{lexicosum} Let $(P,\preceq)$ be a \emph{finite} poset and for each $x\in P$ let $(I_x,A_x)$ be an $I^*$-algebra. We set $I=\bigoplus_{x\in P}I_x$, $A=\bigoplus_{x\in P}A_x$, and we write elements of $A$ in the form $(a_x)_{x\in P}$. We define 
$$\bigL_{x\in P}I_x=\{a\in I|\forall x,y\in P, x\prec y\Rightarrow \max\sigma(a_x)\le \min\sigma(a_y)\}$$
Then $\bigL_{x\in P}I_x$ is an isocone of $A$.
\end{theorem}
\begin{demo}
For ease of notation let $J=\bigL_{x\in P}I_x$, and for a self-adjoint element $a\in A$, let us write $\max(a)$ and $\min(a)$ instead of $\max\sigma(a)$ and $\min\sigma(a)$.

It is clear that $J$ contains $\RR.1_A$ and is norm-closed. Let $a,b\in J$. For every $x,y\in P$ such that $x\prec y$, we have

\be
\max(a_x+b_x)\le\max(a_x)+\max(b_x)
\ee
and 
\be
\min(a_y)+\min(b_y)\le \min(a_y+b_y)
\ee
Using these two inequalities and the definition of $a\in J$ and $b\in J$ we obtain $a+b\in J$.

Consider now some $f\in I(\RR)$. If $x,y\in P$ are such that $x\prec y$, we have :

\bea
\max(\sigma(f(a)_x))&=&\max(\sigma(f(a_x))), \mbox{ obviously}\cr
&=&\max(f(\sigma(a_x))),\mbox{ by the spectral mapping theorem}\cr
&=&f(\max(\sigma(a_x))),\mbox{ since }f\mbox{ is non-decreasing}\cr
&\le& f(\min(\sigma(a_y))),\mbox{ since }a\in J\mbox{ and }f\mbox{ is non-decreasing}\cr
&\le&\min(\sigma(f(a)_y)),\mbox{ by the same steps as above}
\eea
Finally let $a\in I$. We are going to show that $a\in J-J$. For this let $m=\inf_{x\in P}\min(a_x)$. Then $a'=a-m\ge 0$. Let $F : P\rightarrow \RR$ be the function defined by 

$$F(x)=\sum_{z\prec x}\max(a_z')$$

Since $a'\ge 0$, $F$ is an isotone function. Let $f\in J$ and $b\in I$ be defined by $f_x=F(x)1_x$ for all $x\in P$, where $1_x$ is the unit of $A_x$, and $b=a'+f$.

Let $x,y\in P$ be such that $x\prec y$. Then

\bea
\min(b_y)&=&\min(a_y')+F(y)\cr
&=&\min(a_y')+F(x)+\max(a_x')+\sum_{z\prec y, z\not\preceq x}\max(a_z')\cr
&=&\min(a_y')+\max(b_x)+\sum_{z\prec y, z\not\preceq x}\max(a_z')\nonumber
\eea

Since $a'\ge 0$, we have $\min(b_y)-\max(b_x)\ge 0$. Hence, $b\in J$.

Thus $a=a'+m=b+m-f\in J-J$. Hence $I\subset J-J$, thus $I-I\subset J-J$, and consequently $\reel(A)\subset \overline{J-J}$. 
\end{demo}

\begin{propo} With the notations of theorem \ref{lexicosum}, we have
$$P(A)\approx \bigL_{x\in P}P(A_x)$$
as toposets, where $P(A)$ is equipped with the ordering $\preceq_J$, with $J=\bigL_{x\in P}I_x$, and $P(A_x)$ is equipped with $\preceq_{I_x}$.
\end{propo}
\begin{demo}
It is standard that $\coprod_{x\in P}P(A_x)$ is homeomorphic to $P(A)$ when we map the pure state $\phi\in P(A_x)$ to $\tilde \phi$ such that $\tilde\phi=\phi$ on $A_x$ and $\tilde\phi=0$ on $A_y$, $y\not=x$. From now on we identify the two spaces.

Let us show that the ordering is the claimed lexicographic ordering. Consider $x,y\in P$ and   $\phi\in P(A_x)$,  $\psi\in P(A_y)$.

If $x\prec y$, for any $a=(a_z)_{z\in P}\in J$ we have $\phi(a)=\phi(a_x)\le\max(\sigma(a_x))$ and $\psi(a)=\psi(a_y)\ge\min(\sigma(a_y))$. Since $a\in J$ we have $\phi(a)\le \psi(a)$. Hence $\phi\preceq_J\psi$.

Now suppose $x=y$. Then $\phi\preceq_J\psi\Leftrightarrow \forall a\in J$, $\phi(a_x)\le\psi(a_x)$. In order to prove that it is equivalent to $\phi\preceq_{I_x}\psi$, we just need to prove that the projection $\pi_x : A\rightarrow A_x$ is such that $\pi_x(J)=I_x$. For this take $b\in I_x$ and consider $a\in A$ defined by $a_y=\min(\sigma(b))$ if $y\prec x$, $a_y=\max(\sigma(b))$ if $x\prec y$, $a_x=b$ and $a_y=0$ if $y\parallel x$. Clearly $a\in J$ and $\pi_x(a)=b$.

Finally suppose $x\parallel y$. Then choose $b\in I_x$ such that $\sigma(b)\subset]0;+\infty[$ and use the construction above. We obtain an $a\in J$ such that $a_x=b$ and $a_y=0$. Hence $\phi(a)>0$ and $\psi(a)=0$. Since we can exchange the roles of $\phi$ and $\psi$, we see that $\phi\parallel\psi$.
\end{demo}

%It is easy to prove along the same lines that if $a\in J$ is non-derogatory,  $(\sigma(a),\preceq_I)$ is the lexicographic sum over $P$ of the toposets $(\sigma(a_x),\preceq_{I_x})$. (When $a$ has multiple eigenvalues, the ordering is obtained from the previous one by collapsing the corresponding elements.) 

Using the ordering $=$ on $P$, the construction above gives us access to finite direct sums of $I^*$-algebras. What about infinite sums ? Let $S$ be a set of indices and $A_s$ be a unital $C^*$-algebra for any $s\in S$. Then  ${\cal A}:=\{(a_s)_{s\in S}|\sup\|a_s\|<\infty\}$ is a $C^*$-algebra. Let $I_s$ be an isocone in $A_s$, and ${\cal I}=\{(a_s)_{s\in S}\in {\cal A}|\forall s\in S,\ a_s\in I\}$. Then ${\cal I}$ is clearly a pre-isocone. Do we have $\overline{I-I}=A$ ?

\begin{figure}
\begin{center}
\includegraphics[scale=0.7]{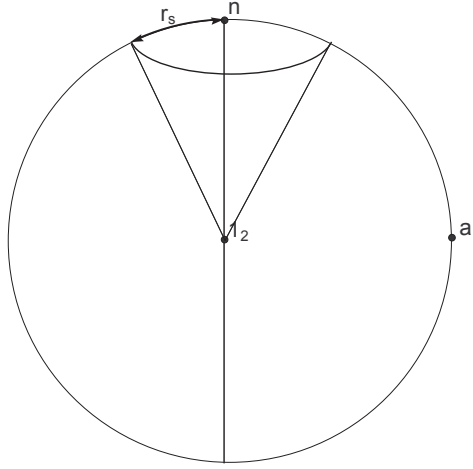}
\end{center}
\caption{The ``closing-cones phenomenon''. The radius $r_s\rightarrow 0$ as $s\rightarrow \infty$.}\label{ccpheno}
\end{figure}

Let us consider the case where $S=\NN$, $A_s=M_2(\CC)$ for all $s$, and $I_s$ is the isocone generated by a spherical cap of $S$ of radius $r_s$ centered on a fixed point $n$, and suppose $r_s\rightarrow 0$ as $s\rightarrow \infty$. Take $a_s=a$ for all $s$, with $a$ a fixed hermitian matrix of trace $1$ which does not lie on the line joining $1_2/2$ and $n$. Suppose for definiteness that $a\in S$ and makes a right angle with $n$ and $1_2/2$ (see figure \ref{ccpheno}).

Then for any $s$, whatever $u_s,v_s\in I_s$ we take such that $u_s-v_s=a$, we have $\|u_s\|\rightarrow \infty$. The same thing happens if we just require that $\|u_s-v_s-a\|<\epsilon$ for some fixed $\epsilon$. Hence $\overline{I-I}\not=\reel(A)$. We need some hypothesis to forbid this ``closing-cone phenomenon'' to happen.

\begin{propo} Let us suppose that there exists $r>0$ such that $\forall s$, $I_s$ contains a ball of radius $r$ and center $c_s$ such that $\sup_s\|c_s\|<\infty$ (non-closing hypothesis). Then, using the above notations,  ${\cal I}$ is an isocone of ${\cal A}$.
\end{propo}  

We leave the easy proof of this proposition to the reader. We note that using only the algebra $M_2(\CC)$ in the sum, we can use it to give some non-trivial examples of infinite-dimensional noncommutative $I^*$-algebras. The next proposition will provide  a different kind of example.

\begin{propo}\label{projinv} Let $A,B$ be $C^*$-algebras, $\pi : A\rightarrow B$ a $*$-morphism, $J$ a pre-isocone  of $B$, and define $I=\reel (\pi^{-1}(J))$. Then :
\begin{enumerate}
\item $I$ is a  pre-isocone of $A$,
\item\label{c2} if $\pi$ is surjective and $J$ is an isocone then $I$ is an isocone.
\end{enumerate}
\end{propo}
\begin{demo}
The first claim is essentially trivial. For the second one, we need to prove that $\overline{I-I}=\reel(A)$. Take $a\in \reel(A)$. Since $\pi(a)\in \reel(B)=\overline{J-J}$ and $\pi$ is surjective, find $a_n,b_n\in I$ such that $\|\pi(a-a_n+b_n)\|<\epsilon$. Now $\pi$ is surjective so that $B\simeq A/\ker\pi$. Thus $$\|\pi(a-a_n+b_n)\|=\inf_{k\in\ker\pi}\|a-a_n+b_n+k\|<\epsilon$$
Hence there exists $k_n\in\ker\pi$ such that $\|a-a_n+b_n+k_n\|<\epsilon$. Moreover, $k_n$ can be supposed to be self-adjoint. Indeed, write $k_n=k+ik'$, with $k$ the self-adjoint and $ik'$ the anti-self-adjoint parts of $k_n$. Then $k$ and $k'$ each belong to $\ker\pi$ and $\|a-a_n+b_n+k\|\le \|a-a_n+b_n+k+ik'\|$ since $\|x\|\le\|x+iy\|$ for any self-adjoint elements $x,y$ in a $C^*$-algebra. Now $\reel(\ker\pi)\subset I$ thus $b_n+k_n\in I$. This proves that $a$ can be arbitrarily approximated by elements of $I-I$.
\end{demo}

Here are two examples where we can use this proposition.

\begin{enumerate}
\item Let $X$ be a compact set and let $A={\cal C}(X,M_2(\CC))$ be the $C^*$-algebra of continuous functions from $X$ to $M_2(\CC)$. Then for any $x\in X$, the evaluation map $\epsilon_x : A\rightarrow M_2(\CC)$ is surjective, which allows us to pull any isocone $I_K$ of $M_2(\CC)$ back to $A$.
\item The Toeplitz algebra  comes equipped with a surjective morphism $\pi$ onto $S^1$. Any toposet structure on $S^1$ then gives rise to a noncommutative $I^*$-algebra structure on the Toeplitz algebra.
\end{enumerate}

\begin{propo}\label{prop6} With the same notations as above, if $J$ is an isocone and $\pi$ is surjective then as toposets we have
$$(P(A),\preceq_I)\simeq (P(B),\preceq_J)+(P^K(A),=)$$
where $K=\ker\pi$, and $P^K(A)$ is the set of pure states which do not vanish on $K$.
\end{propo}
\begin{demo}
We know (see for instance \cite{dix} p.63) that $P(A)=\pi^*P(B)\coprod P^K(A)$ and that  $\pi^* : \phi\mapsto \phi\circ \pi$ is a bijection from $P(B)$ to $\pi^*P(B)$.

Now let $\omega\in P_K(A)$ and $\phi\in P(A)$. Suppose $\omega\perp\phi$   and take $k\in K$. Since $\pm k\in K\subset I$, we deduce that $\omega(k)=\phi(k)$.

Now if $\phi\in \pi^*P(B)$, then $\phi(K)=0$ which entails $\omega(K)=0$, which is absurd. Thus $P_K(A)\parallel\pi^* P(B)$.

If $\phi\in P^K(A)$ then $\phi$ and $\omega$ coincide on $K$, whence they are equal.

Finally if $\phi,\phi'\in \pi^*P(B)$, then $\phi=\psi\circ \pi$. The same goes for $\phi'$, and for any $a\in I$
$$\phi(a)\le\phi'(a)\Leftrightarrow \psi(\pi(a))\le\psi'(\pi(a))\Leftrightarrow \psi\preceq_J\psi'$$
since $\pi$ is surjective.
\end{demo}

The above proposition shows that the construction \ref{c2} in proposition \ref{projinv} when applied to commutative algebras is dual to the embedding of a compact toposet $M$ into a compact set $N$, extending the order on $M$ by the trivial ordering (equality) on $N$.  Suppose now that we have a compact toposet $N$ and a closed subset $M\subset N$. The restriction of functions gives rise to a projection morphism $\pi : {\cal C}(N)\rightarrow {\cal C}(M)$ such that $\pi(I(N))$ is exactly the isocone $I(M,\preceq_M)$ where $\preceq_M$ is the restriction of the ordering on $N$ to $M$. That $\pi(I(N))$ is closed is a consequence of the existence of a ``Tietze extension theorem'' in the category of toposets (see \cite{bes1}). In the noncommutative setting, we can prove the following proposition.

\begin{propo}\label{projdir}
Let $A,B$ be $C^*$-algebras, $\pi : A\rightarrow B$  a $*$-morphism, and $I$ a pre-isocone  of $A$ . Then :
\begin{enumerate}
\item $\overline{\pi(I)}$ is a pre-isocone  of $B$. Moreover, if $A$ is finite-dimensional, then $\pi(I)$ is closed. 
\item If $I$ is an isocone and $\pi$ is surjective, then $\overline{\pi(I)}$ is an isocone  of $B$.
\end{enumerate}

\end{propo}
\begin{demo}
Clearly $\overline{\pi(I)}$ is a closed convex cone containing the constants. Take $f\in I(\RR)$ and $y\in \overline{\pi(I)}$. Then $y=\lim_n\pi(a_n)$, $a_n\in I$. We have

\bea
f(y)&=&f(\lim_n\pi(a_n))\cr
&=&\lim_n f(\pi(a_n)),\mbox{ by continuity of functional calculus}\cr
&=&\lim_n \pi(f(a_n)),\mbox{ since }\pi\mbox{ is a }*\mbox{-morphism}\cr
&\in&\overline{\pi(I)},\mbox{ since }I\mbox{ is a pre-isocone}
\eea
%In the strong case, $y=\lim_n\pi(a_n)$, $z=\lim_n\pi(b_n)$, $y\wedge z={1\over 2}(\lim_n(\pi(a_n+b_n)-|\lim_n(\pi(a_n-b_n))|=\lim_n\pi(a_n\wedge b_n)$, by continuity arguments.

To prove that $\pi(I)$ is closed in the finite-dimensional case, take $y$ in the closure of $\pi(I)$ and consider a sequence $\pi(x_n)$ which converges towards $y$, with $x_n\in I$. Take an increasing homeomorphism from $\RR$ onto $]-1;1[$. Then $f(x_n)\in I\cap  B$, where $B$ is the unit ball of $A$. Since $A$ is finite-dimensional, we can consider a subsequence $x_{n_k}$ such that $f(x_{n_k})$ converges to some element $z$. Since $f(x_{n_k})\in I$ and $I$ is closed, $z\in I$. Hence $\pi(f(x_{n_k}))=f(\pi(x_{n_k}))$ converges to $\pi(z)\in \pi(I)$. Thus $f(y)=\pi(z)\in \pi(I)$. Now $\pi(I)$ is stable by isotone functional calculus, hence $y=f^{-1}(f(y))$ belongs to $\pi(I)$ (the fact that $f^{-1}$ is not defined on $\RR$ is not a problem since it is defined on $\sigma(f(y))$).
\end{demo}

It is immediate by definition that $\pi^*P(B)$ equipped with the restriction of $\preceq_I$ is isomorphic to $P(B)$ with the ordering defined by $\overline{\pi(I)}$.

\underline{Remark}: we do not know any example where $\pi(I)$ is not closed.

Finally, we note that by using successively propositions \ref{projdir} and \ref{projinv} we get an interesting corollary.

\begin{cor} Let $A$ be a sub-$C^*$-algebra of $B$,  $I$ an isocone of $A$ and $K$ a closed two-sided ideal of $B$. Then $\overline{I+\reel(K)}$ is an isocone of the $C^*$-algebra $A+K$. 
\end{cor}
\begin{demo} That $A+K$ is a $C^*$-algebra is a classical result (\cite{kad}, p 717). For the rest we note that $I$ is a pre-isocone of $B$, hence $\overline{\pi(I)}$ is a pre-isocone of $B/K$, with $\pi : B\rightarrow B/K$ the canonical surjection. Thus $\reel(\pi^{-1}(\overline{\pi(I)}))$, which is easily seen to be equal to $\overline{I+\reel(K)}$, is a pre-isocone of $B$. Moreover from $\overline{I-I}=\reel(A)$ we get $\overline{I-I+\reel(K)}=\reel(A+K)$.
\end{demo}

If we view the elements of $I$ as noncommutative isotonies and those of $K$ as negligible in some way, this corollary tells us that we can perturb an isotony with something negligible and still get an isotony. To make this a little bit more explicit, consider the following example. We take a compact toposet $M$, and let $H=L^2(M)$ for a regular Borel measure $\mu$ on $M$, % wikipedia page on continuous functions on compact hausdorff spaces. En fait il suffit que la fonction caractéristique de l'ouvert où une fonction continue ne s'annule pas soit mesurable. Dans ce cas pour tout f non nulle on peut trouver Psi=1_cet ouvert tq f\Psi non nulle, donc i(f) n'est pas nul. Il suffit donc que les ouverts soient mesurables.
then the represention by pointwise multiplication $\iota : M\rightarrow {\cal B}(H)$, $f\mapsto (\psi\mapsto f\psi)$ is faithful, so that we can identify ${\cal C}(M)$ with $A:=\iota({\cal C}(M))$ and $I(M)$ with $I:=\iota(I(M))$. Taking for $K$ the ideal of compact operators ${\cal K}(H)$ in $B={\cal B}(H)$, we see that the set of isotonies on $M$ perturbed with compact operators is still an isocone (note that in this case one can prove that $I+\reel(K)$ is closed). This was certainly to be expected since compact operators play the role of infinitesimals in noncommutative geometry. As a final remark, let us observe what proposition \ref{prop6} becomes in this case. The pure states which vanish on $K$ are the evaluation maps on $A$, that is Dirac delta functions. On the other hand, those which don't are vector states of the form $a\mapsto\langle \psi,a\psi\rangle$ for some $\psi\in H$ of unit $L^2$ norm. In physics language the elements of $P(A)$  are wave functions on $M$. Those of the Dirac type are completely localized, and these are the only ones which are comparable to each other for the order relation induced by $I+\reel (K)$.
 
\section{Finite-dimensional $I^*$-algebras}
\subsection{Overview}
Consider a finite-dimensional $C^*$-algebra $A$. We know that $A$ is isomorphic to a direct sum of full matrix algebras, and we will always consider such an isomorphism as given. That is to say, we consider $A$ to be of the form

$$A=\bigoplus_{x=1}^kM_{n_x}(\CC)\subset M_N(\CC),\quad N=\sum_{x=1}^kn_x$$

Suppose that $I$ is an isocone of $A$. Then we are going to prove that $I$ is of the form $\bigL_{x\in P}I_x$, where $P$ is the set $\{1;\ldots;k\}$ equipped with a poset structure. Note that if $n_1=\ldots=n_k=1$, $A$ is $\CC^k$ and $I=I(P)$. Hence we can loosely describe this result as saying that every finite noncommutative poset is a poset each point of which has been given an internal ordered structure. Moreover, we will see that these internal structures can be enumerated : if $n_x\not=2$, $I_x=\Herm(n_x)$ is the only possibility, and if $n_x=2$, $I_x=I_{K_x}$ where $K_x$ is some compact and geodesically convex subset of $S^2$ with non-empty interior.

In what follows we consider an isocone $I\subset\reel(A)\subset\Herm(N)$ and we denote by $\interior(I)$ the interior of $I$ for the topology of $A$ (which must be non empty).

Here is a sketch of the different stages through which we will arrive at the classification theorem.

In subsection \ref{innerordering} we show that $I$ defines an ordering on the spectrum of its interior elements, and that this ordering is locally constant. In fact it is truly constant when $A=M_N(\CC)$, and in that case we call it the inner ordering of $I$. The crucial property is that the inner ordering is trivial if and only if $I$ is trivial.

In subsection \ref{projinisoc} we prove that the projections in $I$ form a lattice (for the usual meet and join for projections), and this will be our main tool in proving the theorem. We will also prove, using projections, that the Hasse diagram of the inner ordering of an isocone of $M_N(\CC)$ cannot be disconnected, unless it is trivial.

In subsection \ref{caseN3} we prove that the isocones of matrix algebras $M_N(\CC)$ with $N\ge 3$ are all trivial. We do this by induction on $N$, the difficult part being the case $N=3$. In this case, we show that the inner ordering of $I$ is trivial by exhibiting a particular element in $I$, constructed thanks to the lattice property of the projections, such that the ordering on its spectrum has a disconnected Hasse diagram. This will prove the triviality of $I$.

Finally in subsection \ref{classif} we prove that every finite-dimensional isocone is a lexicographic sum.
%We will come to the theorem in several stages. 
% For any $r\in\NN$ we call $H_r$ the affine hyperplane of matrices $a$ such that $\tr(a)=r$. It is orthogonal to the line $\RR.1_N$.

\subsection{The inner ordering of an isocone}\label{innerordering}
Here is a heuristic discussion of what we are aiming at in this subsection. Take an element $a\in I$. Recall that $I\cap C^*(a)$ is an isocone of $C^*(a)\simeq \C(\sigma(a))$. This endows $\sigma(a)$ with a toposet structure. Basically the idea is that if there is no sudden change of dimension of $C^*(a)$ when $a$ is moved around a little bit in $I$, then the toposet structure must be constant. Of course  to make sense of this idea we need to have some means of identifying the different spectra $\sigma(a)$ when $a$ varies.

For $a\in S_N$, call $s_a : \dcg 1..N\dcd\rightarrow \sigma(a)$ the map such that $s_a(1)<\ldots<s_a(N)$. Then $s_a^*$ is an isomorphism between $\C(\sigma(a))$ and $\C(\dcg 1..N\dcd)$, which we hereafter identify with $\CC^N$. The morphism $s_a^*$ has all imaginable properties (for us it will be a $*$-isomorphism, hence an isometry). By composition with the Gelfand-Naimark isomorphism $\theta_a : C^*(a) \rightarrow \C(\sigma(a))$ we have a $*$-isomorphism $\phi_a=\theta_a\circ s_a^* : C^*(a)\rightarrow \CC^N$. This isomorphism takes the following simple matricial form : let $U$ be a unitary matrix such that $U^*aU=\diag(a_1,\ldots,a_N)$ with $a_1<\ldots<a_N$, then $\phi_a(x)=U^*xU$, where we identify a diagonal matrix with an element of $\CC^N$.

The image of $I\cap C^*(a)$ by $\phi_a$ is an isocone of $\CC^N$ that we call $I(a)$. It is associated with a partial order on $\dcg 1..N\dcd$ which we denote by $\le_a$.

%, by $I_1^+$ the intersection of $I$ with $\ge 0$ matrices of trace $1$.

%Similarly we denote by $I(a)_1^+$ the elements of $I(a)$ which are $\ge 0$ and sum to $1$ (that is the intersection of $I$ with the standard $N-1$ simplex of $\RR^N$).

%We observe that $I=\RR_+I_1^++\RR\id$ and $I(a)=\RR_+I(a)_1^++\RR.1$.

We will need two lemmas. The first can be expressed by saying that if we take a fixed element $f$ in $\CC^N$, the map $a'\mapsto \phi_{a'}^{-1}(f)$ gives an identification of $f$ with an element of $C^*(a')$ which is continuous around $a$ if $a$ is non-derogatory.

\begin{figure}[hbtp]
\begin{center}
\includegraphics[scale=0.5]{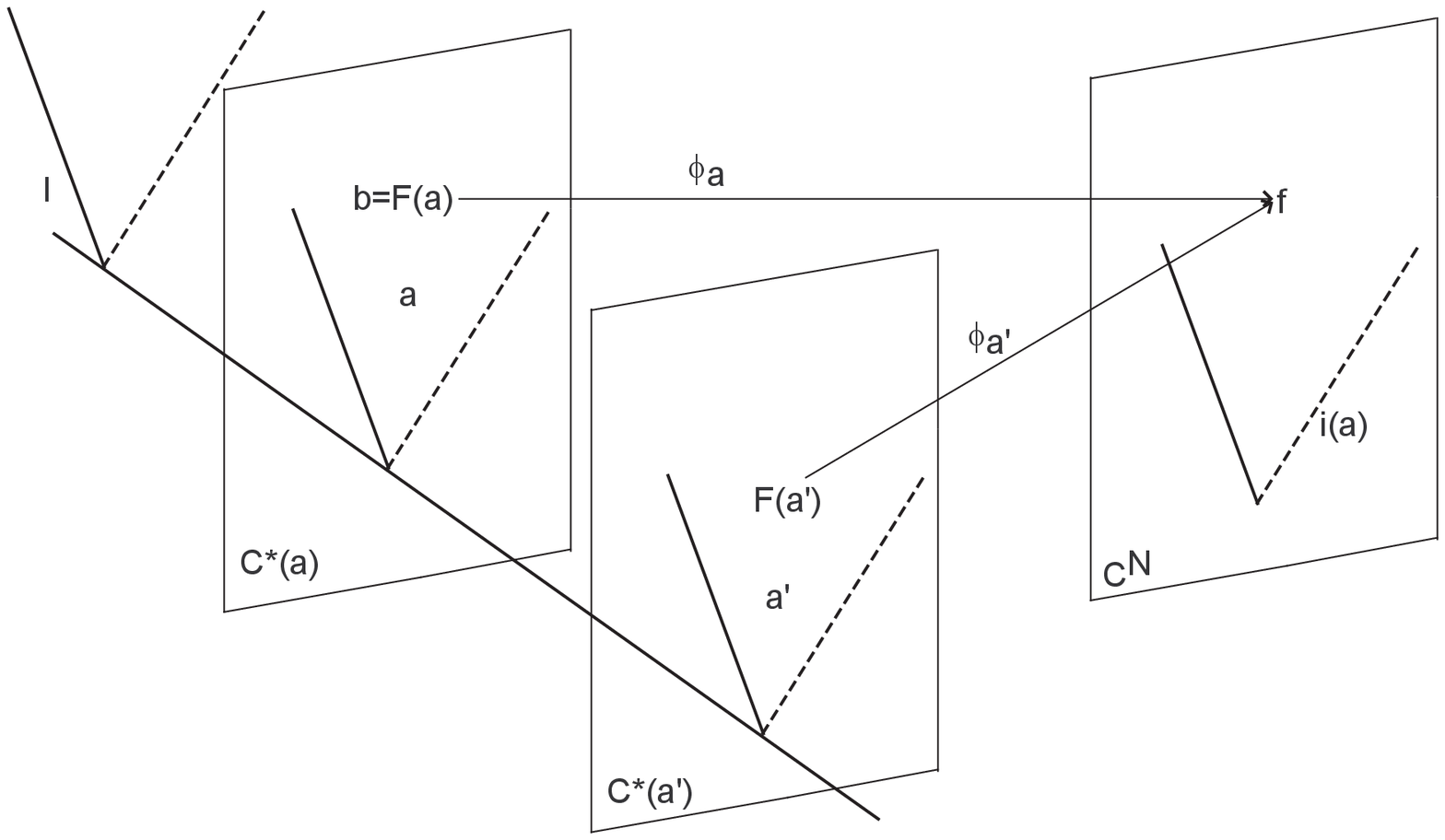}
\end{center}
\end{figure}

\begin{lemma} (continuous identification lemma) Let $a\in S_N$ and let $b\in C^*(a)$. Then there exists a continuous function $F\in\C(\RR,\RR)$ such that $F(a)=b$ and an $\epsilon>0$ such that $\|a-a'\|<\epsilon\Rightarrow \phi_{a'}(F(a'))=\phi_a(F(a))$.
\end{lemma}
\begin{demo}
Les $\lambda_1<\ldots<\lambda_N$ be the eigenvalues of $a$, $e_1,\ldots,e_N$ be the corresponding eigenbasis, and $b_i$, $1\le i\le n$ be the eigenvalue of $b$ corresponding to $e_i$. Let $\delta$ be the infimum of the distance between two eigenvalues of $a$. Let $F$ be a continuous function on $\RR$ such that $F(]\lambda_j-\delta/3;\lambda_j+\delta/3[)=b_j$. Then for $\epsilon>0$ small enough, if $\|a-a'\|<\epsilon$, the $j$-th largest eigenvalue $\lambda_j'$ of $a'$ will be  in the interval $]\lambda_j-\delta/3;\lambda_j+\delta/3[$, hence we will have (with obvious notations) $a=\sum_i \lambda_i p_{e_i}$, and $a'=\sum_i \lambda_i'p_{e_i'}$, $F(a)=\sum_i b_i p_{e_i}$, and $F(a')=\sum_i b_i p_{e_i'}$. Thus $\phi_{a'}(F(a'))=\phi_a(F(a))=\phi_a(b)=(b_1,\ldots,b_N)$.
%Alternate proof : the same result can presumably be achieved by writing $a=U_aD_aU^*_a$, $b=U_aD_b U^*_a$ and $F: U_{a'}D_{a'}U_{a'}^*\mapsto U_{a'} D_bU_{a'}^*$ which is certainly continuous for $a'$ close enough to $a$. 
\end{demo}

The second lemma is about convex sets. Its proof is left to the reader.

\begin{lemma}\label{lemmaconv} Let $I$ be a convex set in some normed vector space $V$. Let $a,b\in I$ with $a\in\interior(I)$, and let $W$ be a vector (or affine) subspace such that $a,b\in W$. If $b$ is in the  interior of $W\cap I$ with respect to the relative topology of $W$, then $b\in\interior(I)$.
\end{lemma}
% \begin{demo} Since $b$ is the relative interior of $W\cap I$, there exists an open ball $B(b,\epsilon)$ such that $B(b,\epsilon)\cap W\subset I$. The intersection of the line $(ab)$ with $B(b,\epsilon)$ is a segment $]b'b''[$, which lies in $W\cap B(b,\epsilon)$, hence in $I$. Now $a$ is in the interior of $I$, so consider an open ball $B(a,\delta)\subset I$. Since $I$ is convex, the convex hull of $B(a,\delta)\cup]b'b''[$ is  entirely in $I$, is open (it is the interior of a finite cone with apex $b'$ or $b''$), and contains $b$. Hence $b\in\interior(I)$.

% Alternate proof under the assumption that $I$ is closed : suppose $b\in\partial I$. The set $B(b,\epsilon)\cap W$ of the previous proof being in $I$, it shows that $W$ must be a supporting subspace of $I$. But this is contradictory with $a$ being in the interior of $I$.
% \end{demo} 

We can now state and prove the theorem.

\begin{theorem}\label{ouvert} For every $a\in\interior(I)\cap S_N$, there exists an $\epsilon>0$ such that $\forall a' \in A$, $\|a-a'\|<\epsilon\Rightarrow$ $a'\in \interior(I)\cap S_N$ and $\le_{a'}=\le_a$.
\end{theorem}
\begin{demo} The first part of the statement is trivial since $\interior(I)\cap S_N$ is open, but it is needed to write it in order for the second part to make sense.

We will prove that for $a'$ sufficiently close to $a$, $I(a)=I(a')$.

For this take $f\in\interior(I(a))$ and call $b\in I\cap C^*(a)$ its pre-image under $\phi_a$. The second lemma tells us that  $b\in\interior(I)$.   

Let $F$ be a continuous real functions such as in lemma 1, and a corresponding $\epsilon$. Take an open ball $B(b,\delta)\subset I$. By continuity of the functional calculus, there is an open ball $B(a,\epsilon')$ such that $F(B(a,\epsilon'))\subset B(b,\delta)$. Hence for $a'$ such that $\|a-a'\|<\epsilon''=\min(\epsilon,\epsilon')$ we have $F(a')\in B(b,\delta)\subset I$, thus $F(a')\in C^*(a')\cap I$, and $\phi_{a'}(F(a'))=f$, which proves that $f$ stays in $I(a')$ for $a'$ sufficiently close to $a$. Hence we have shown that $\interior(I(a))\subset I(a')$, which in our case obviously entails that $I(a)\subset I(a')$.

Conversely, take $f\notin I(a)$ and let $b$ be such that $\phi_a(b)=f$. Then $b\notin I$. Consider a function $F$ given by the continuous identification lemma. Since $I$ is closed, there is an open ball $B(b,\delta)$ disjoint from $I$. By continuity of the functional calculus, its pre-image contains an open ball $B(a,\epsilon)$. With $\epsilon$ small enough we have $\forall a'\in B(a,\epsilon)$, $\phi_{a'}(F(a'))=f$. Moreover, we have $F(a')\notin I$, hence $f\notin I(a')$. Thus for $a'$ close enough to $a$, we have $\RR^N\setminus I(a)\subset \RR^N\setminus I(a')$.
\end{demo}

Now consider the set $P_N$ of all possible orderings on $\dcg 1..N\dcd$ which are at most as fine as the natural ordering. It is a finite set partially ordered by inclusion. For a given ordering $R\in P_N$, let $U_R(I)$ be the set $\{a\in \interior(I)\cap S_N|\le_a=R\}$. The theorem above shows that $U_R(I)$ is open. Hence the open set $\interior(I)\cap S_N$ is a finite union of disjoint open sets, which must be unions of connected components of $\interior(I)\cap S_N$. Let $C_N$ be the complement of $\S_N$ inside $\interior(I)$, that is, the set of derogatory matrices in $I$. If $A=M_N(\CC)$, then $C_N$ is a manifold of codimension 3, hence $\interior(I)\cap S_N=\interior(I)\setminus C_N$ is connected.

% Indeed, let $\phi$ be a diffeomorphism between $\interior(I)$ and $\RR^{N^2}$. Then $\phi(\interior(I)\setminus C_N)=\RR^{N^2}\setminus \phi(C_N)$ is connected. It only remains to show that $\interior(I)$ is diffeomorphic to $\RR^{N^2}$. But this is obvious : every convex open set is diffeomorphic to ambiant space, and this is what the interior of $I$ is.

We therefore obtain the following theorem (and definition):

\begin{theorem}
Let $I$ be an isocone in $M_N(\CC)$. For all $a\in\interior(I)\cap S_N$, the ordering $\le_a$ is constant. We call it the \emph{inner ordering} defined by $I$.
\end{theorem}

% \underline{Remark 2} : using remark 1 and this theorem, we see that  the inner ordering of $I$ is total iff no projection is in the interior of $I$ as soon as $N\ge 3$.

% Here is a direct elementary proof of this claim.

% \begin{lemma} $\interior(I)\cap S_N$ is connected.
% \end{lemma}
% \begin{demo} We show first that $U=\reel M_N(\CC)\cap S_N$ is connected. Take $a\in U$ and let $d$ be the diagonal matrix containing the eigenvalues of $a$ in increasing order. We have $a=udu^*$ with $u\in U(N)$, and since $U(N)$ is path-connected, we have a path from $t\mapsto u(t)du(t)^*$ from $a$ to $d$. Now take $b\in U$. We can similarly connect $b$ a diagonal matrix $d'$ with increasing eigenvalues. It only remains to show that there is a path from $d$ to $d'$, but it is obvious since the set of increasing sequences in $\RR^N$ is  convex.

% Now take an homeomorphism $\phi : \interior(I)\rightarrow \RR^{N^2}$. 
% \end{demo}

Let us look at the case of isocones of the form $I_K$ in $M_2(\CC)$ to see what can happen :

\begin{itemize}
%\item Consider first the commutative case. In that case the components of $\interior(I)\cap S_N$  are separated by the hyperplanes of equation $x_i=x_j$. Hence they are intersections of half-spaces, and as such are convex. Each element $a$ of $\interior(I)\cap S_N=I\cap S_N$ generates $\RR^N$, hence $\le_a=\le_I$ is constant. No point of the boundary belongs to $\S_N$.
\item Suppose $K$ contains no antipodal points. Here the order is constant everywhere except on the scalar matrices, where the spectrum degenerates. %Note that here $\interior(I)\cap S_N$ is not the same as $I\capS_N$, since there are points in the boundary with two different eigenvalues. Moreover, there are some rank one projectors in $\interior(I)$, showing that we can points of $C_N$ inside $\interior(I)$.
\item Suppose now that $K$ is a hemisphere. In that case every element of the interior of $I_K$ has simple eigenvalues. We then have $U_{1\oplus 2}(I_K)=\interior(I)$. If $a\in \partial I$, then either $\sigma(a)$ degenerates to a single point, or, if $a$ is not a multiple of the identity, it remains a pair of points but the ordering degenerates to the trivial one. \emph{This shows that a degeneration of the ordering can happen at the boundary of $I$, even if the eigenvalues stay simple.}
\end{itemize}

With these examples we see that $\le_a$ can be or not be constant on $I\cap S_N$. So the theorem is optimal in a way. However, we can notice at once that the converse part of the proof of theorem \ref{ouvert} does not use the fact that $a\in\interior(I)$. Hence we see that $\le_{a}$ can only denegerate by becoming less fine that the constant ordering defined by the elements of $\interior(I)\cap S_N$. We also see that $\le_a$ is constant on every connected components of $\partial I\cap S_N$.

We now come to a simple but important consequence of the previous theorem.

\begin{propo}\label{triv} Let $I$ be an isocone in $M_N(\CC)$. The inner ordering is trivial iff $I$ is trivial.
\end{propo}
\begin{demo} The `if' part is obvious. If the inner ordering is trivial, then for all $a\in I\cap S_N$, $-a\in I$ since the function $x\mapsto -x$ is isotone with respect to $\le_a$. Hence the vector space $I\cap (-I)$ has non empty interior. Thus $I\cap (-I)=I=-I=\Herm(N)$.
\end{demo}

\subsection{Projections in isocones}\label{projinisoc}

Projections play a particularly important role in finite-dimensional isocones thanks to proposition \ref{combproj}. Let us introduce some notations.

If $a,b,c,\ldots$ are any vectors or subsets of a vector space, we will write $[a,b,c,\ldots]$ for the vector subspace generated by $a,b,c,\ldots$. If $V$ is a vector subspace of $\CC^N$, $p_V$ will be the orthogonal projection on $V$, with the exception that we write $p_x$ instead of $p_{[x]}$ when $x$ is a single non-zero vector.

We write ${P}_k(N)$ for the set of rank $k$ projections in $\Herm(N)$ and ${  P}_k(I)$ those which are in $I$.

We consider $p_L$ and $p_N$ two projections in an isocone $I\subset A$, and  we first look for an eigendecomposition of a convex combination $tp_L+(1-t)p_N$, $t\in]0;1[$. Then, using the isotone functional calculus, we will see that many projections in $I$ can be found starting with $p_L$ and $q_L$ (16 in general).

The main tool is Halmos' two subspaces theorem \cite{halmos2}, which we recall here in the form given in \cite{BS}.
  
\begin{theorem}\label{2subspaces} (Halmos) Let $L$ and $N$ be two closed subspaces of a Hilbert space $\H$. Let us write

$$L=(L\cap N)\buildrel\perp\over\oplus (L\cap N^\perp)\buildrel\perp\over\oplus L_0$$
$$L^\perp=(L^\perp\cap N)\buildrel\perp\over\oplus (L^\perp\cap N^\perp)\buildrel\perp\over\oplus L_0'$$
 
If one of the spaces $L_0,L_0'$ is non trivial (i.e. the spaces $L,N$ are in general position), then both these spaces have the same dimension and there exists a unitary operator $R : L_0'\rightarrow L_0$  and a hermitian $H : L_0\rightarrow L_0$ satisfying $0\le H\le \id_{L_0}$, with $\ker H=\ker (\id_{L_0}-H)=0$, such that

$$p_L=\id_{L\cap N}\oplus \id_{L\cap N^\perp}\oplus 0_{L^\perp\cap N}\oplus 0_{L^\perp\cap N^\perp}\oplus U^*  \pmatrix{\id_{L_0}&0\cr 0&0} U$$
and
$$p_N=\id_{L\cap N}\oplus 0_{L\cap N^\perp}\oplus \id_{L^\perp\cap N}\oplus 0_{L^\perp\cap N^\perp}\oplus U^*  \pmatrix{\id_{L_0}-H&W\cr W&H} U$$

with $U=\diag(\id_{L_0},R)$,  and $W=(H-H^2)^{1/2}$. Note that the matrix in the middle of the product is an endomorphism of $L_0\oplus L_0$. The product itself is an endomorphism of $L_0\buildrel\perp\over\oplus L_0'$.
\end{theorem}

Let us make two observations.
\begin{enumerate}
\item This theorem is valid in infinite dimension. In finite dimension, the part about the kernel of $H$ implies that $\sigma(H)\subset]0;1[$, which is not necessarily true in infinite dimension.
\item The hypothesis that $L$ and $N$ are in general position, in other words that $L_0$ and $L_0'$ are non trivial, is exactly equivalent to the hypothesis that $p_L$ and $p_N$ do not commute (see \cite{BS} prop 1.5).
\end{enumerate}
 
Let us now take $t\in]0;1[$ and look for the spectral decomposition of the operator $tp_L+(1-t)p_N$.  Since the computation is elementary (and can be found in \cite{BS} for the case $t=1/2$), we only summarize it.

We have :

\bea
tp_L+(1-t)p_N&=& \id_{L\cap N}\oplus t\id_{L\cap N^\perp}\oplus (1-t)\id_{L^\perp\cap N}\oplus 0_{L^\perp\cap N^\perp}\cr
& &\oplus U^*  \pmatrix{\id_{L_0}-(1-t)H&(1-t)W\cr (1-t)W&(1-t)H} U\nonumber
\eea

Writing $T=\pmatrix{\id_{L_0}-(1-t)H&(1-t)W\cr (1-t)W&(1-t)H}-{1\over 2}\pmatrix{\id_{L_0}&0\cr 0&\id_{L_0}}$ we see easily that the spectrum of $T$ is of the form $K_t\cup(-K_t)$, where $K_t$ is a finite subset of $]0;1/2[$.

Writing $S_t^+=K_t+{1\over 2}$ and $S_t^-=-K_t+{1\over 2}$, and $V_t^+,V_t^-$ for the corresponding eigenspaces, the spectrum of $tp_L+(1-t)p_N$ is thus of the form (if some of the subspaces in the decomposition are trivial, just delete the corresponding eigenvalues) :

\begin{enumerate}
\item If $t<1/2$ : $0<S_t^-<t<1-t<S_t^+<1$, corresponding to $L^\perp\cap N^\perp\oplus V_t^-\oplus L\cap N^\perp\oplus L^\perp\cap N\oplus V_t^+\oplus L\cap N$,
\item If $t>1/2$ : $0<S_t^-<1-t<t<S_t^+<1$, corresponding to $L^\perp\cap N^\perp\oplus V_t^-\oplus L^\perp\cap N\oplus L\cap N^\perp\oplus V_t^+\oplus L\cap N$,
\item If $t=1/2$ : $0<S_t^-<1/2<S_t^+<1$, corresponding to $L^\perp\cap N^\perp\oplus V_t^-\oplus (L^\perp\cap N\oplus L\cap N^\perp)\oplus V_t^+\oplus L\cap N$,
\end{enumerate}

% We note that $V_t^+,V_t^-,L_0,L_0'$ all have the same dimension. Moreover,  $V_t^+\cap L_0=V_t^-\cap L_0=V_t^+\cap L_0'=V_t^-\cap L_0'=0$. %Hence $L_0\oplus L_0'=L_0\oplus V_t^+$.  

Now we start using isotone functional calculus to obtain new projections in $I$ from the decomposition above. We could continue working with projections, but for ease of notations let us introduce $H_I$, the set of subspaces such that the associated projections belongs to $I$.  Here is the list of spaces that belong to $H_I$ thanks to isotone functional calculus :

\begin{enumerate}
\item\label{enu1} $L\cap N$,
\item\label{enu2} $(L\cap N)\oplus V_t^+$, for any $t$
\item\label{enu3} $(L\cap N)\oplus V_t^+\oplus (L\cap N^\perp)$, for $t>1/2$
\item\label{enu4} $(L\cap N)\oplus V_t^+\oplus (L\cap N^\perp)\oplus (L^\perp\cap N)$, for any $t$
\item\label{enu5} $(L\cap N)\oplus (L\cap N^\perp)\oplus (L^\perp\cap N)\oplus L_0\oplus L_0'=L+N$ (since $V_t^++V_t^-=L_0+L_0'$),
\item\label{enu6} $(L\cap N)\oplus V_t^+\oplus  (L^\perp\cap N)$, for $t<1/2$.
\end{enumerate}
And of course $H_I$ also contains
\begin{enumerate}\setcounter{enumi}{6}
\item\label{enu7} $L$,
\item\label{enu8} $N$.
\end{enumerate}

The first consequence is the following important property :

\begin{propo}\label{propo9}
Let ${\mathcal L}(\H)$ be the lattice of subspaces of $\H$. Then $H_I$ is a sublattice of ${\mathcal L}(\H)$ (hence ${ P}(I)$ is a sublattice of ${P}(\H)$).
\end{propo}
\begin{demo}
Let $p_L,p_N$ belong to $I$. If they commute, then $p_L\vee p_N$ and $p_L\wedge p_N$ belong to $I$ by the weak isocone property. Hence $L\cap N$ and $L+N\in H_I$.

If they do not commute, use \ref{enu1} and \ref{enu5}.
\end{demo}

We will now use this result to make various combinations of subspaces in the preceding list using $\cap$ and $+$. It will be useful to introduce the decomposition obtained by exchanging the roles of $L$ and $N$. The notations are summarized in table \ref{tablesub}.

We note the trivial but useful fact that $L_0+L_0'=N_0+N_0'$ and we call $W$ this common space. 

Note also that $L_0=L\cap W$ and $N_0=N\cap W$ entail that $L_0\cap N_0=(L\cap N)\cap W=0$. Moreover, the four spaces $L_0,N_0,L_0',N_0'$ having the same dimension, we have $L_0\oplus N_0=W$.

Similarly, $N_0\cap L_0'=(N\cap W)\cap (L^\perp\cap W)=(N\cap L^\perp)\cap W=0$, and $N_0+L_0'=W$, and also $N_0'\oplus L_0=W$, $N_0'\oplus L_0'=W$.

\begin{table}[hbtp]
\begin{center}
\begin{tabular}{|c|c|c|}
\hline
$L\cap N$ & $L\cap N^\perp$ & $L_0$\\
\hline
$L^\perp\cap N$& $L^\perp\cap N^\perp$ & $L_0'$\\
\hline
\end{tabular} \hspace{1cm}\begin{tabular}{|c|c|c|}
\hline
$L\cap N$ & $L^\perp\cap N$ & $N_0$\\
\hline
$L\cap N^\perp$& $L^\perp\cap N^\perp$ & $N_0'$\\
\hline
\end{tabular}
\caption[smallcaption]{ Two decompositions of $H$ into orthogonal summands}
\label{tablesub}
\end{center}
\end{table}

\begin{theorem}
Let $p_L$ and $p_N$ be two non commuting projections in $I$. Then the sublattice generated by $L\cap N+L\cap N^\perp$, $L\cap N+L^\perp\cap N$, $L\cap N+L_0$, $L\cap N+N_0$ belongs to $H_I$. This sublattice is distributive, and is isomorphic to the lattice of subsets of $\{L\cap N^\perp ;L^\perp\cap N;L_0;N_0\}$. 
\end{theorem}
\begin{demo}
For the sake of simplicity of notations, let $O=L\cap N$, $A_1=L\cap N^\perp$, $A_2=L^\perp\cap N$, $A_3=L_0$ and $A_4=N_0$. Let $\LL$ be the lattice generated by $O+A_i$, $i=1..4$, and $\LL'$ be the lattice consisting of the corresponding projections.

The set $\{p_O,p_{A_i}|i=1,..4\}$ is a Foullis-Holland set, that is to say a non-empty subset of an orthomodular lattice such that whenever one chooses three distinct elements of this set, one of them commutes with the other two. Such a set generates a distributive lattice $\LL''$ (see \cite{greechie}). In a distributive lattice, all elements can be reduced to the normal form $p_{i_1}\vee\ldots \vee p_{i_k}$. Now our lattice $\LL'$ is the sublattice of $\LL''$ consisting of elements larger than $p_O$. Thus we just have to check that there are exactly 16 elements in this lattice, that is to say that no two subpaces of the form $O+\sum_i A_i$ are equal. This is easy by direct inspection.

Now the only thing that remains to be proven is that $\LL\subset H_I$. For this it suffices to prove that the generators belong to $H_I$. Using the list of subspaces just before proposition \ref{propo9} we obtain :

\begin{itemize}
\item the intersection of type \ref{enu3} and type \ref{enu7} is  $L\cap ((L\cap N)\oplus V_t^+\oplus (L\cap N^\perp))=(L\cap N)\oplus (L\cap N^\perp)$. %Indeed, $L\cap ((L\cap N)\oplus V_t^+\oplus (L\cap N^\perp))\subset L\cap (L\oplus L_0')=L$.

\item Similarly (exchanging $L$ and $N$ or working with type \ref{enu6} and \ref{enu8}) we find that $(L\cap N)\oplus (L^\perp\cap N)\in H_I$.

\item  \ref{enu3} + \ref{enu7} gives $(L\cap N)+(L\cap N^\perp)+V_t^++L_0=(L\cap N)\oplus (L\cap N^\perp)\oplus L_0\oplus L_0'=L\oplus L_0'\in H_I$ (since $V_t^++L_0=L_0+L_1$). Then \ref{enu6} + \ref{enu8} gives $N\oplus N_0'\in H_I$.  Now $(L\buildrel\perp\over\oplus L_0')\cap (N\buildrel\perp\over\oplus N_0')=((L\cap N)+W)\buildrel\perp\over\oplus (L\cap N^\perp))\cap ((L\cap N)+W)\buildrel\perp\over\oplus (L^\perp\cap N))=L\cap N+W$.   Intersecting $L\cap N+W$ and $L$ we get $(L\cap N)\buildrel\perp\over\oplus L_0\in H_I$.
\item Symmetrically we obtain $(L\cap N)\buildrel\perp\over\oplus L_0'\in H_I$.
\end{itemize}
\end{demo}
 
The following is an easy corollary.

\begin{propo}\label{cororder}
For all $t\in]0;1[$, the order induced by $I$ on $\sigma(tp_L+(1-t)p_N)$ is at most as fine as the one given by the following Hasse diagram :

$$\xymatrix@R-15pt{
 & 1 & \\
 & S_t^+\ar@{-}[u]& \\
\ t\ \ar@{-}[uur]&  &1-t\ar@{-}[uul]\\
 & S_t^-\ar@{-}[uu]&\\
 & 0 \ar@{-}[uul]\ar@{-}[u]\ar@{-}[uur]&
}
$$
\end{propo}
\begin{demo}
We know from the previous theorem that $p_{L\cap N+L\cap N^\perp}$ and $p_{L\cap N+L^\perp\cap N}$ both belong to $I$. They also belong to $C^*(tp_L+(1-t)p_N)$, hence $1-t$ and $t$ are not comparable and $1-t,t$ are not $\le S_t^+$. Similarly, using the space $L\cap N+W$ we see that $S_t^-$ cannot be $\le 1-t,t$.
\end{demo}

We close this section with a topological property of the set of projections in an isocone $I$ of $M_N(\CC)$.

\begin{lemma}\label{localsurj} 
Let $a\in\interior(I)\cap S_N$, let $h(a)$ be an eigenprojection of rank $k$ of $a$ such that $h(a)\in I$. Then there exists $\epsilon>0$ such that $P_k(N)\cap B(h(a),\epsilon)\subset I$.
\end{lemma}
\begin{demo}
Let $a_1<\ldots<a_N$ be the eigenvalues of $a$. For ease of notation, we can suppose without loss of generality that $h(a)$ is the eigenprojection corresponding to the $k$ first eigenvalues of $a$. We begin by extending the definition of $h$ : it is clear that there exists a continuous real function $\tilde h$, such that $\tilde h(a_1)=\ldots=\tilde h(a_k)=1$ and $\tilde h(a_i)=0$ for $i>k$, and such that for all $j$, $a_j$ lies in an open interval on which $\tilde h$ is constant. Then $\tilde h(a)=h(a)$, and $\tilde h$ is defined for any hermitian matrix.

Now for any $b\in M_N(\CC)$ define the conjugation map $c_b :  U(N)\rightarrow M_N(\CC)$, $U\mapsto U^*bU$.

We observe that for $\eta>0$ small enough, $\tilde h$ defines a mapping from $B(a,\eta)$ to $P_k(N)$, and $B(a,\eta)\subset \interior(I)\cap S_N$. Moreover, we can supppose thanks to theorem \ref{ouvert} that $\tilde h(B(a,\eta))\subset I$. Finally, since $\tilde h(x)$ is insensitive to the spectrum of $x$, as long as $x$ is close enough to $a$, one has for $\eta$ small enough, $\tilde h(B(a,\eta))=\tilde h(B(a,\eta)\cap {\mathcal O}(a))$, where ${\mathcal O}(a)=c_a(U(N))$ is the unitary orbit of $a$. 
 
By continuity of $c_a$, there is an open set $W$ containing $1_N$ such that $c_a(W)=B(a,\eta)\cap {\mathcal O}(a)$. Now since $\tilde h\circ c_a=c_{h(a)}$ on $W$, it is sufficient to show  that $c_{h(a)}$ is an open mapping from $U(N)$ onto  ${\mathcal O}(h(a))$ for the relative topology of  ${\mathcal O}(h(a))$. But $c_{h(a)}$ is a continuous surjective map from $U(N)$ onto ${\mathcal O}(h(a))$, which defines a continuous bijection $\tilde c$ from $U(N)/{\rm Stab}(h(a))$ to ${\mathcal O}(h(a))$. Since $U(N)$ is compact, so is $U(N)/{\rm Stab}(h(a))$, hence $\tilde c$ is an homeomorphism. Since the quotient is one of topological groups, the quotient map is open, hence $c_{h(a)}$ is open.
% Let $\epsilon>0$ be given. Then there exists $\eta>0$ such that for all unitary $U$,  $\|U-I\|<\eta\Rightarrow \|U^*aU-a\|<\epsilon$. 
\end{demo}

This property will be used crucially in what follows.

\begin{cor}\label{dis} Let $I$ be a non-trivial isocone of $M_N(\CC)$ and $\preceq$ be its inner ordering. Then the Hasse diagram of $\preceq$ is connected.
\end{cor}
\begin{demo}
Suppose that the Hasse diagram is not connected. Take a component $C$. Then both the functions $\delta_C$ and $-\delta_C$ are isotone for $\preceq$. Consequently, for any $a\in\interior(I)\cap S_N$, there will be some projector $h(a)\in I$, for which we also have $-h(a)$ in $I$. Now we know from the previous lemma that set of all such $h(a)$ contains a bit of a unitary orbit, hence it will contain a matrix basis $(m_i)$. Thus $m_i$ and $-m_i$ belong to $I$, and $I$ is a convex cone, hence $I=\Herm(N)$.
\end{demo}

\subsection{Isocones in $M_N(\CC)$ for $N\ge 3$}\label{caseN3}
This section consists of a single theorem.
\begin{theorem}
Let $I$ be an isocone  of $M_N(\CC)$ with $N\ge 3$. Then $I=\Herm(N)$.
\end{theorem}
\begin{demo}
We first consider the case $N=3$.

Take $a\in\interior(I)\cap S_3$. Set $a=U\diag(a_1,a_2,a_3)U^*$, with $a_1<a_2<a_3$, $p_3=U\diag(0,0,1)U^*$ and $p_{23}=U\diag(0,1,1)U^*$. We know that $p_3$ and $p_{23}\in I$. Fix a $t\in]0,1[$, for instance $t=1/2$. Then $b:=tp_3+(1-t)p_{23}\in \interior(I)$ (from lemma \ref{lemmaconv}). We are going to show that we can find an element  $b'$ which is still in $\interior(I)$ and is still a convex combination of a rank one and a rank two projection, but with respective ranges in general position. We will then conclude by proposition \ref{cororder} and corollary \ref{dis}.

Fix an $\epsilon>0$ such that $B(b,\epsilon)\subset I$. Take a non zero $\theta\in\RR$ and let $V=U\pmatrix{\cos\theta&0&-\sin\theta\cr 0&1&0\cr \sin\theta&0&\cos\theta}$ and $p_3'=V\diag(0,0,1)V^*$. If $u_1,u_2,u_3$ are the columns of $U$, then $p_3$ is the projection on the line $[u_3]$ and $p_3'$ is the projection on the line generated by $v_3=-(\sin\theta) u_1+(\cos\theta) u_3$. By taking $\theta$ sufficiently small we can assure that 

%By taking $V$ sufficiently close to $U$ in $U(3)$, and $p_3'=V^*\diag(0,0,1)V$, we can assure that 

%Now let $u_1,u_2,u_3$ be the columns of $U$ and  $\theta\in\RR$. Let $V=U\pmatrix{\cos\theta&0&-\sin\theta\cr 0&1&0\cr \sin\theta&0&\cos\theta}$ and $p_3'=V\diag(0,0,1)V^*$. By taking $\theta$ sufficiently small we can assure that 

\begin{itemize}
\item $\|p_3-p_3'\|<\epsilon/t$,
\item $p_3'\in I$,
\item $p_3'$ does not commute with $p_{23}$.
\end{itemize}

The first property is obvious and the second one follows at once from lemma \ref{localsurj}. The third one is an easy calculation. This calculation provides a formal proof of the intuitively obvious fact that for $\theta$ small    $[v_3]=$ range of $p_3'$ and $[u_2,u_3]=$ range of $p_{23}$ are in general position.

We now write $L$ for the range of $p_3'$ and $N$ for the one of $p_{23}$ in order to use the notations of theorem \ref{2subspaces} and table \ref{tablesub}. We have $L=L_0$ since $L$ is one-dimensional, and $N_0'=N^\perp$ also for dimensional reason. Therefore  $L\cap N=L\cap N^\perp=L^\perp\cap N^\perp=0$ and $L_0\simeq L_0'$ and $L^\perp\cap N$ all have dimension $1$.

Now we set $b'=tp_3'+(1-t)p_{23}$. We have $\|b'-b\|=t\|p_3-p_3'\|<\epsilon$. Thus $b'\in\interior(I)$. Moreover, using proposition \ref{cororder}, we see that the corresponding ordering on $\sigma(b')$ is at most as fine as $(1\oplus 3)+2$, where by $1,2,3$ we mean the eigenvalues of $b'$ in ascending order (namely the single element of  $S_t^-$, $1-t$, and the single element of $S_t^+$).  Consequently, it has a disconnected Hasse diagram.

Thus $I$ is trivial by corollary \ref{dis}.

% Now write $L$ for the range of $p_3'$ and $N$ for the one of $p_{23}$. Since $p_{23}$ and $p_3'$ do not commute we know that $L$ and $N$ are in general position. Then, using the notations of theorem \ref{2subspaces} and table \ref{tablesub}, we have $L=L_0$ since $L$ is one-dimensional, and $N_0'=N^\perp$ also for dimensional reason. Therefore  $L\cap N=L\cap N^\perp=L^\perp\cap N^\perp=0$ and $L_0\simeq L_0'$ and $L^\perp\cap N$ all have dimension $1$.
% The dimension box of $b'$ is 

% \begin{table}[h]
% \begin{center}
% \begin{tabular}{c|c|c}
% 0 & 0 & 1 \\
% \hline
% 1 & 0 & 1 \\
% \end{tabular}
% \caption[dimension box]{dimension box of a convex combination of a rank one and a rank two projection in $M_3(\CC)$.}
% \label{boxn3p1q2}
% \end{center}
% \end{table}

We consider now the case $N=4$ and pick an $a\in\interior(I)\cap S_4$. Up to replacing $I$ with $UIU^*$, for some unitary matrix $U$, we can as well suppose that $a=\diag(a_1,\ldots,a_4)$ with $a_1<\ldots<a_4$. We call $B$ the $C^*$-subalgebra $M_3(\CC)\oplus \CC\subset M_4(\CC)$. Since $B$ contains  $a\in\interior(I)$, the pre-isocone $I\cap B$   has a non-empty interior in $B$, and thus is an isocone of $B$.

Let $\pi : B\rightarrow M_3(\CC)$ be the first projection.  By   proposition \ref{projdir} above, $J=\pi(I\cap B)$ is an isocone of $M_3(\CC)$, hence it is trivial.

Consequently, $J$ contains $\diag(1,0,0)$, $\diag(0,1,0)$ and $\diag(0,0,1)$. Thus $I$ contains elements of the form $\diag(1,0,0,c)$ and $\diag(0,1,0,c')$ for some $c$ and $c'$, which both act as isotone functions for the order $\preceq_a$. Hence $1\parallel 2$ for $\preceq_a$. Similarly we can show that $1\parallel 3$ and $2\parallel 3$. But we can do exactly the same reasoning using $B'=\CC\oplus M_3(\CC)$ instead of $B$. We thus come to the conclusion that $2,3,4$ are also incomparable with each other. Hence the Hasse diagram of the inner ordering is disconnected, $2$ being incomparable with every other element. This shows that $I$ is trivial.

The proof goes on by induction.
\end{demo}

\subsection{The classification theorem}\label{classif}

We use the following notations : $\pi_x : \bigoplus_{i\in P}M_{n_i}(\CC)\rightarrow M_{n_x}(\CC)$ is the projection on the $x$-th summand, $\pi_{x,y} :  \bigoplus_{i\in P}M_{n_i}(\CC)\rightarrow M_{n_x}(\CC)\oplus M_{n_y}(\CC)$ is $\pi_x\oplus\pi_y$ and $I_x=\pi_x(I)$. Let $N=n_1+\ldots+n_k$. %The  unit of $M_N(\CC)$ is written $1_N$, and its zero $0_N$. Moreover $\id_x$ will denote the element of $A$ which has $1_{n_x}$ on the $x$-th summand and zero everywhere else. 

\begin{theorem}
Let $I$ be an isocone in the finite-dimensional $C^*$-algebra $A=\bigoplus_{x\in P}M_{n_x}(\CC)$, with $P=\{1;\ldots;k\}$, $k\in\NN^*$, $n_x\in\NN^*$. Then there exists a poset structure on $P$ such that $I=\bigL_{x\in P} I_x$. 
\end{theorem}

The rest of this section is devoted to the proof of this theorem. We begin with the $k=2$ case, that is where $A=M_n(\CC)\oplus M_p(\CC)$.

\begin{lemma} Let $a=(a_1,a_2)\in\interior(I)\subset M_n(\CC)\oplus M_p(\CC)$ and let $\lambda_1=\max(\sigma(a_1))$ and $\lambda_2=\max(\sigma(a_2))$.
\begin{itemize}
\item If $\max(\sigma(a))=\lambda_1$ and $\lambda_1$ has multiplicity one, then $\iota_1=1_n\oplus 0_p\in I$.
\item If $\max(\sigma(a))=\lambda_2$ and $\lambda_2$ has multiplicity one, then $\iota_2=0_n\oplus 1_p\in I$.
\end{itemize}
\end{lemma}
\begin{demo} Let us prove the first claim. Call $e_\lambda\oplus 0$ the eigenvector of $\CC^n\oplus\CC^p$ corresponding to the eigenvalue $\lambda$ of $a$, and $p_{e_\lambda}$ the corresponding rank one projector. Since the condition of the lemma must hold on a neighbourhood of $a$, we get a familly $e_{\lambda_i}$, $i=1,\ldots,n$ of generating vectors for $\CC^n$. Since $p_{{e_\lambda}_i}\oplus 0_p\in I$ for all $i$, we have by proposition \ref{propo9} that  $p_{[e_{\lambda_1},\ldots,e_{\lambda_n}]}\oplus 0_p=\iota_1\in I$.
\end{demo}

We observe that at least one of the two cases in the lemma has to occur for at least one element of $\interior(I)$. Now we claim that 

\begin{itemize}
\item If $\iota_1$ and $\iota_2$ are in $I$, then $I=I_1\oplus I_2$.
\item If $\iota_1\notin I$ and $\iota_2\in I$, then $I=\bigL_{x\in 1\oplus 2}I_x$.
\item If $\iota_2\notin I$ and $\iota_1\in I$, then $I=\bigL_{x\in 2\oplus 1}I_x$.
\end{itemize}

Suppose $\iota_1,\iota_2\in I$. The inclusion $I\subset I_1\oplus I_2$ being alway true, we prove the converse. Let $a_1\in I_1$ and $a_2\in I_2$. Then by definition there exist $b_1,b_2$ such that $a_1\oplus b_2\in I$ and $b_1\oplus a_2\in I$. Let $\lambda$ be a constant such that $(a_1+\lambda)\oplus (b_2+\lambda)$ and $(b_1+\lambda)\oplus (a_2+\lambda)$ are both positive. Then since $\iota_1$ and $\iota_2$ commute with  $(a_1+\lambda)\oplus (b_2+\lambda)$ and $(b_1+\lambda)\oplus (a_2+\lambda)$, their products with these elements are in $I$. Hence $(a_1+\lambda)\oplus 0$ and $0\oplus(a_2+\lambda)$ are in $I$, thus their sum is. We can then substract $\lambda 1_N$.

Next we suppose that $\iota_1\notin I$. We suppose moreover that there exists $a=a_1\oplus a_2$ in $I$ such that $\max(a_1)>\min(a_2)$, and we intend to arrive at a contradiction. We remark first that we can take $x\in\interior(I)\cap S_N$, and replace $a$ with $(1-\epsilon)a+\epsilon x$, $\epsilon>0$ in order to displace all possible equalities among eigenvalues. Hence, we can suppose without loss of generality that $a$ is non-derogatory and belongs to the interior of $I$. By hypothesis we have $\iota_1\notin I$, hence $\max(a_1)<\max(a_2)$. Let us call $p_a$ the eigenprojection  corresponding to $\max(a_1)$ and $V(a)$ the eigenspace of $a_2$ corresponding to the all the eigenvalues $\ge \max(a_1)$. By isotone calculus, we have $p_a\oplus p_{V(a)}\in I$. Using conjugation by unitaries of the form $U_1\oplus 1_p$, we can move $a_1$ around while keeping $a_2$ constant. Hence we obtain $p_i\oplus p_{V(a)}\in I$ for a family of projections $p_i$ corresponding to a basis of $\CC^n$, and using the lattice property of ${ P}(I)$, we conclude that $1_n\oplus p_{V(a)}\in I$. Now this is true for all $a'$ in a neighbourhood of $a$. Moving $a$ this time with unitaries of the form $1_n\oplus U_2$, we can obtain a family of projections of the form $1_n\oplus p_{V_i}$, with the subspaces $V_i$ in general position, and taking the infimum of these elements, we conclude that $1_n\oplus 0_p\in I$.

The third case is of course symmetric to the second one. 

Now that we are finished with the $k=2$ case, we will use it to prove the general result. The proof is a variation on the Kakutani-Stone theorem, which is combinatorially more complex because of the noncommutativity, but topologically simpler because of the finite dimensionality.

Let $\le$ be defined on $P$ by $x\le y$ iff $x=y$ or $\forall a\in I$, $\max(a_x)\le \min(a_y)$. This relation is  antisymmetric : if $\forall a\in I$, $\max(a_x)\le \min(a_y)\le\max(a_y)\le\min(a_x)$ we have $a_x=\lambda 1_{n_x}$ and $a_y=\lambda 1_{n_y}$ for all $a\in I$, which is impossible since $I$ has a non-empty interior. Since $\le$ is obviously transitive and reflexive, it is a partial order on $P$.

\begin{lemma} Let $x\in P$ and $f^x$ be the function on $P$ defined by $f^x(y)=1$ if $x\le y$, $f^x(y)=0$ otherwise. Then $(f^x(y))_{y\in P}\in I$.
\end{lemma}
\begin{demo} Suppose $u\in P$ is such that $x\not\le u$. Then by the $k=2$ case, $1_{n_x}\oplus 0_{n_u}\in \pi_{x,u}(I)$. In other words, there exists $a^u\in I$ such that $a^u_x=1_{n_x}$, $a^u_u=0_{n_u}$.

On the other hand, if $x\le u$, we can find $a^u\in I$ such that $a^u_x=1_{n_x}$ and $a^u_u=1_{n_u}$ (for instance $a^u=1_N$).

Let $H\in I(\RR)$ be such that $H(]-\infty;0])=0$, $H([0;1])=[0;1]$ and $H([1;+\infty[)=1$. Up to a composition with $H$, we can suppose without loss of generality that $a^u$ is such that $0\le a^u_y\le 1_{n_y}$ for all $y\in P$, and $a^u_y=1_{n_y}$ if $x\le y$.

Now let $a$ be the average of all $a^u$ for $u\in P\setminus\{x\}$. If $x\not\le z$ then $a_x$ has its spectrum inside $[0;{n-1\over n}]$. If $x\le z$, we have obviously $a_z=1_{n_z}$. We get the needed element of $I$ by composing with an affine non-decreasing function which vanishes on $]-\infty;{n-1\over n}]$ and takes the value $1$ on $[1;+\infty[$.
\end{demo}

Now let $L:=\bigL_{x\in P}I_x$, we obviously have $I\subset L$, we must prove the converse. Thanks to proposition \ref{combproj} we only need to do it  for projections.

\begin{lemma}
Let $p=(p_x)_{x\in P}$ be a projection in $L$. For all $x,y\in P$, there exists $a^{xy}\in I$ such that : $a^{xy}_x=p_x$, $a^{xy}_y=p_y$, $a^{xy}_z=0_{n_z}$ or $1_{n_z}$ for every $z\not=x,y$.
\end{lemma}
\begin{demo} First step : we find an element $b\in I$ such that $b_x=p_x$ and $b_y=p_y$ (possible thanks to the $k=2$ case).

Second step : $c=H(b)$ where $H$ is the same function as in the previous lemma. We note that whenever $p_x\not=0$ and $x< z$, then $c_z=1_{n_z}$, and whenever $p_x\not=1_{n_x}$ and $z<x$, then $c_z=0_{n_z}$ (the same goes for $y$).

Third step : we use the lemma above to find an $f^x$ and an $f^y$. If $p_x$ and $p_y$ are non zero we set $d=f^x\vee f^y$ which belongs to $I$ since $f^x$ and $f^y$ commute. If $p_x=0$ and $p_y\not=0$, we set $d=f^y$, and symmetrically if $p_x\not=0$ and $p_y=0$ (the case $p_x=p_y=0$ is trivial).

Last step : We call $a^{xy}=dc$. We note that $d$ and $c$ commute and are $\ge 0$. Hence $dc\in I$.

We easily check that $a^{xy}_x=p_x$ and $a^{xy}_y=p_y$. Let us check the other condition.

First case : suppose $p_x$ and $p_y$ do not vanish. Then if $x<z$ we have $c_z=d_z=1_{n_z}$.  The same goes if $y<z$. If neither $x\le z$ nor $y\le z$, then $d_z=0$, hence $a^{xy}_z=0$.

Second case : $p_x$ does not vanish, $p_y=0$. Then if $x<z$ we have $a^{xy}_z=1_{n_z}$ as above, and if $x\not\le z$, we have $d_z=0$. 
\end{demo}

We can conclude the proof of the theorem : let $p$ be a projection in $L$. Then for all $x,y\in P$ we find an $a^{xy}$ as above. Now we set

\begin{itemize}
\item $a^x=\sup\{a^{xy}|y\in P\setminus\{x\}\}$,
\item $a=\inf\{a^x|x\in P\}$
\end{itemize}

We note that all the elements of which we take the supremum in the first formula commute with each other and belong to $I$. Hence so does $a^x$. Moreover $a^x_x=p_x$ and either $a^x_y=1_{n_y}$ or $a^x_y=p_y$ for $y\in P$. Hence the $a^x$ commute with each other (because for a particular $y$ the different $a^x_y$ are either equal to a constant or a constant operator). Hence $a\in I$. Finally since $a_x$ is the infimum of a familly of operators one of them being equal to $p_x$ and the other being equal to $p_x$ or $1_{n_x}$, we have $a_x=p_x$. Thus $a=p$.

\section{Conclusion and outlook}

The classification theorem might seem to be a bit disappointing at first, since it shows that there are not that many interesting examples in finite dimension. To this, several answers can be given. Firstly, finite-dimensional $C^*$-algebras are not particularly interesting either, but this does not account for the richness of the full theory. Our result can thus be viewed as an invitation to explore the infinite dimensional case, starting with   almost commutative algebras. In fact, a classification result in this case is already within reach \cite{BB}. It would also be interesting to know whether the causal cones of \cite{francoeckstein} are stable under isotonies, in which case they would be isocones. Moreover, the very fact that noncommutative examples are hard to find but nonetheless exist is rather encouraging: it indicates that our set of axioms is just constraining enough to be consistent but not trivial.

From another point of view, the finite-dimensional case does not merit to be so easily dismissed. We have already emphasized its importance in Connes' approach to unification as well as in   causal set quantum gravity. The introduction of a noncommutative inner structure for the points of a causal set would be a most natural step linking the two theories. A natural question is  whether the causal structure would then necessarily remain purely commutative. Our classification result seems to give an interesting mixture of affirmative and negative answer, and this with very few physical input (no dynamics). The non-trivial part of the partial order is seen to be mostly commutative, the only exception coming from the $M_2(\CC)$ summands. Even though it is nothing more that a wild guess at this stage, it is hard not to notice the similarity with the strong CP problem.

Finally let us mention that the context of complex $C^*$-algebras  is possibly too restrictive. Already in the Chamseddine-Connes spectral model we see real $C^*$-algebras appearing, and this pleads for an extension of our studies in this direction. Another context which   deserves some attention is the one of Jordan algebras. In this setting it is indeed possible to find examples which go beyond our classification: there exists non-trivial isocones  in the so-called Jordan algebras of Clifford type, also known as spin factors, even in infinite dimension. A classification result here is also in progress.

\end{document}